\documentclass[12pt]{article}
\usepackage{amsmath}
\usepackage{times}
\usepackage{graphicx}
\usepackage{color}
\usepackage{multirow}
\usepackage[authoryear]{natbib}
\usepackage{rotating}
\usepackage{bbm}
\usepackage{latexsym}

\usepackage{setspace}
\usepackage[font=doublespacing]{caption}

\usepackage{geometry}
\usepackage[utf8]{inputenc} 
\usepackage[T1]{fontenc}    
\usepackage{hyperref}\hypersetup{
    colorlinks=true,
    linkcolor=black,
    filecolor=blue,      
    urlcolor=black,
    citecolor=black,
}
\usepackage{url}            
\usepackage{booktabs}       
\usepackage{amsfonts}       
\usepackage{nicefrac}       
\usepackage{microtype}      
\usepackage{xcolor}         

\usepackage{colortbl}
\usepackage[flushleft]{threeparttable}
\usepackage{amsmath}\allowdisplaybreaks
\usepackage{amsthm}
\usepackage{algorithm,algorithmic}
\usepackage{multirow}
\usepackage{enumitem,picinpar}
\usepackage{subfigure}
\usepackage{pifont}
\usepackage{array}
\usepackage{graphicx}
\usepackage{amssymb}

\newtheorem{theorem}{Theorem}

\newtheorem{definition}{Definition}
\newtheorem{corollary}{Corollary}
\newtheorem{remark}{Remark}

\newtheorem{assum_abbr}{A}

\newtheorem{innercustomgeneric}{\customgenericname}
\providecommand{\customgenericname}{}
\newcommand{\newcustomtheorem}[2]{%
  \newenvironment{#1}[1]
  {%
   \renewcommand\customgenericname{#2}%
   \renewcommand\theinnercustomgeneric{##1}%
   \innercustomgeneric
  }
  {\endinnercustomgeneric}
}

\newcustomtheorem{custheorem}{Theorem}
\newcustomtheorem{cuslemma}{Lemma}
\newcustomtheorem{cuscorollary}{Corollary}
\newcustomtheorem{cusassum}{Assumption}

\begin{document}


{\bf Obtaining Lower Query Complexities through Lightweight  Zeroth-Order Proximal Gradient Algorithms}

\ \\
{\bf  Bin Gu }\\
{jsgubin@gmail.com}\\
{School of Artificial Intelligence, Jilin University, China}\\
{Mohamed bin Zayed University of Artificial Intelligence,  UAE}\\
{\bf Xiyuan Wei}\\
{xiyuan@tamu.edu}\\
{ Texas A\&M University, USA}\\
{\bf Hualin Zhang}\\
{zhanghualin98@gmail.com}\\
{Mohamed bin Zayed University of Artificial Intelligence,  UAE}\\
{\bf Yi Chang}\\
{yichang@jlu.edu.cn}\\
{School of Artificial Intelligence, Jilin University, China}\\
{\bf Heng Huang}\\
{henghuanghh@gmail.com}\\
{University Maryland College Park, USA}


\thispagestyle{empty}
\markboth{}{NC instructions}
\ \vspace{-0mm}\\
%
\begin{center} {\bf Abstract} \end{center}
Zeroth-order (ZO) optimization  is one key technique for machine learning problems where gradient calculation is expensive or impossible. Several variance reduced ZO proximal  algorithms have been proposed to speed up ZO optimization for non-smooth  problems, and all of them opted for the coordinated ZO  estimator against the random ZO estimator when approximating the true gradient, since the former is more accurate. While the random ZO estimator  introduces bigger error and makes convergence analysis more challenging compared to  coordinated ZO  estimator, it requires only $\mathcal{O}(1)$ computation, which is significantly less than $\mathcal{O}(d)$ computation of the coordinated ZO estimator, with $d$ being dimension of the problem space. To take advantage of the computationally efficient nature of the random ZO estimator, we first propose a ZO objective decrease (ZOOD) property which can incorporate two different types of errors in  the upper bound of  convergence rate. Next, we propose two generic reduction frameworks for ZO optimization which can  automatically derive the convergence results for convex and non-convex problems respectively, as long as the convergence rate for the inner solver satisfies  the ZOOD property. With  the application of two  reduction frameworks on our proposed ZOR-ProxSVRG and ZOR-ProxSAGA, two  variance reduced ZO proximal  algorithms  with fully random ZO estimators, we  improve the state-of-the-art function query complexities  from $\mathcal{O}\left(\min\{\frac{dn^{1/2}}{\epsilon^2}, \frac{d}{\epsilon^3}\}\right)$ to  $\tilde{\mathcal{O}}\left(\frac{n+d}{\epsilon^2}\right)$ under $d > n^{\frac{1}{2}}$ for non-convex problems, and from $\mathcal{O}\left(\frac{d}{\epsilon^2}\right)$ to $\tilde{\mathcal{O}}\left(n\log\frac{1}{\epsilon}+\frac{d}{\epsilon}\right)$ for convex problems. Finally, we conduct experiments to verify the superiority of our proposed methods.

\section{Introduction}

In this paper, we study the following finite-sum optimization problem.
\begin{equation}
\label{problem}
    \begin{aligned}
		\min_{\mathbf{x} \in \mathbb{R}^d} \left\{F(\mathbf{x}) \overset{\textrm{def}}{=} f(\mathbf{x}) + r(\mathbf{x}) \overset{\textrm{def}}{=} \frac{1}{n}\sum_{i=1}^n f_i(\mathbf{x}) + r(\mathbf{x})\right\},
	\end{aligned}
\end{equation}
where $f(\mathbf{x})$ is a smooth function, and $r(\mathbf{x})$ is a relatively simple (but possibly non-differentiable) convex function. Problem \eqref{problem} summarizes  an extensive number of important (regularized) learning problems, e.g., sparse adversarial attack \cite{fan2020sparse}, $\ell_1$-regularized logistic regression \citep{koh2007interior}, robust OSCAR \cite{gu2018inexact}. In many machine learning tasks, such as black-box adversarial attacks on deep neural networks (DNNs) \cite{papernot2017practical,madry2017towards, kurakin2016adversarial}, bandit optimization \cite{flaxman2005online} and reinforcement learning \cite{choromanski2018structured}, calculating the explicit gradients of objective function is expensive or infeasible. As a result, the huge class of first-order gradient algorithms are inapplicable to these tasks. However, zeroth-order (ZO) algorithms use function evaluations to estimate the explicit gradient, which makes them suitable for solving the above tasks. On the other hand, recent research has also  delved into the generalization capabilities of zeroth-order (ZO) algorithms \cite{nikolakakis2022black}. This exploration, combined with classical results regarding convergence guarantees in optimization, renders ZO algorithms both robust and potent for a diverse array of black-box learning problems.

Since zeroth-order algorithms use function evaluations to estimate the explicit gradients to update the solutions, the number of function evaluations for an algorithm to converge is an intuitive metric to assess its performance. Formally, we refer to the number of function evaluations to converge as an algorithm's \textit{function query complexity} (FQC). In this paper, we study zeroth order algorithms with proximal update. In the literature, a proximal stochastic ZO algorithm, RSPGF, was proposed in \cite{ghadimi2016mini}, with FQC of $\mathcal{O}(\frac{d}{\epsilon^4})$ and $\mathcal{O}(\frac{d}{\epsilon^2})$ solving non-convex and convex problems respectively, where $d$ is the dimension of problem. Recently, \cite{kungurtsev2021zeroth} and \cite{chen2018universal} proposed zeroth-order proximal stochastic gradient methods for solving weakly convex stochastic optimization problems, respectively. Further, several variance-reduced variants were proposed for solving non-convex problems. Specifically, ZO-ProxSVRG and ZO-ProxSAGA with FQC of $\mathcal{O}(\frac{d n^{\frac{2}{3}}}{\epsilon^2})$ was proposed in \cite{huang2019faster}, where $n$ is the number of individual functions (\emph{i.e.}, $f_i(\cdot)$ in Eq.\eqref{problem}). PROX-ZO-SPIDER  with  FQC of $\mathcal{O}(\min\{\frac{d n^{\frac{1}{2}}}{\epsilon^2}, \frac{d}{\epsilon^3}\})$ was proposed in \cite{ji2019improved}. Comparison of algorithms is summarized in Tables \ref{t1} and \ref{t2}.
	
\begin{table*}		 	
\caption{Comparison of ZO proximal algorithms solving a non-convex smooth objective with a convex non-smooth regularizer. (VR, CZOG and CR are the abbreviations of variance reduction,   coordinated ZO estimator and convergence rate respectively, $n$ denotes the number of individual functions, $d$ denotes the dimension of problem, $S$ denotes the number of epochs, $m$ denotes the epoch size, $K=mS$, $\mu$ is the smoothing parameter, and $G_*^2 = \max\left\{G_\emph{C}^2, G_\emph{NC}^2\right\}$ with $G_\emph{C}^2$ and $G_\emph{NC}^2$ being constant errors defined in Definition \ref{a6}.)}
	\linespread{1.1}
	\centering
	\begin{small}\setlength{\tabcolsep}{0.9mm}
	\begin{threeparttable}
		\begin{tabular}{lccccc}
			\toprule
			Algorithms & VR & CZOG & CR & System error \tnote{1}   & FQC  \\
			\midrule
			RSPGF \cite{ghadimi2016mini} & No & No & $\mathcal{O}\left(\frac{d^{1/4}}{S^{1/4}}\right)$ & $\mathcal{O}\left( \mu^2\right)$  & $\mathcal{O}\left(\frac{d}{\epsilon^4}\right)$    \\
			ZO-ProxSVRG \cite{huang2019faster}  & Yes  & Yes & $\mathcal{O}\left(\frac{\sqrt{d}}{\sqrt{K}}\right)$& $\mathcal{O}\left( \mu^2\right)$ & $\mathcal{O}\left(\frac{dn^{\frac{2}{3}}}{\epsilon^2}\right)$   \\
			ZO-ProxSAGA \cite{huang2019faster}  & Yes  & Yes & $\mathcal{O}\left(\frac{\sqrt{d}}{\sqrt{S}}\right)$  & $\mathcal{O}\left( \mu^2\right)$  & $\mathcal{O}\left(\frac{dn^{\frac{2}{3}}}{\epsilon^2}\right)$ \\
			PROX-ZO-SPIDER \cite{ji2019improved}  & Yes  & Yes & $\mathcal{O}\left(\frac{1}{\sqrt{K}}\right)$ & $\mathcal{O}\left( \mu^2\right)$ & $\mathcal{O}\left(\min\{\frac{dn^{\frac{1}{2}}}{\epsilon^2}, \frac{d}{\epsilon^3}\}\right)$   \\ \hline
			AdaptRdct-NC (Ours) & & & &  \\
			(ZOR-ProxSVRG/ZOR-ProxSAGA)  & \multirow{-2}{*}{Yes} & \multirow{-2}{*}{No} & \multirow{-2}{*}{$\mathcal{O}\left(\frac{\sqrt{d}}{\sqrt{S}}\right)$}  & \multirow{-2}{*}{$\mathcal{O}\left( \mu^2 + G_*^2\right)$}  &  \multirow{-2}{*}{$\tilde{\mathcal{O}}\left(\frac{n+d}{\epsilon^2}\right)$}  \\
				\bottomrule
			\end{tabular}
			 \begin{tablenotes}
    \item[1] The system error measures the divergence between the lowest functional values  the algorithm converges to and the function itself.
    \end{tablenotes}
    \end{threeparttable}
	\end{small}
	\label{t1}
\end{table*}
	
\begin{table*}
	\caption{Comparison of ZO proximal algorithms solving a convex smooth objective with a convex non-smooth regularizer. ($\rho\in (0, 1)$ denotes a constant,  the other notations are same as Table \ref{t1}.)}
	\linespread{1.1}
	\centering
	\begin{small}\setlength{\tabcolsep}{0.9mm}
		\begin{tabular}{lccccc}
			\toprule
			Algorithms & VR & CZOG & CR  & System error & FQC \\
			\midrule
			RSPGF \cite{ghadimi2016mini} & No & No & $\mathcal{O}\left(\frac{\sqrt{d}}{\sqrt{S}}\right)$  & $\mathcal{O}\left( \mu^2\right)$  & $\mathcal{O}\left(\frac{d}{\epsilon^2}\right)$ \\ \hline
			AdaptRdct-C (Ours) & & & & \\
			(ZOR-ProxSVRG/ZOR-ProxSAGA) & \multirow{-2}{*}{Yes} &  \multirow{-2}{*}{No} & \multirow{-2}{*}{$\mathcal{O}\left(\rho^S\right)$} & \multirow{-2}{*}{$\mathcal{O}\left( \mu^2 + G_*^2\right)$}  &  \multirow{-2}{*}{$\tilde{\mathcal{O}}\left(n\log\frac{1}{\epsilon}+\frac{d}{\epsilon}\right)$}   \\
			\bottomrule
		\end{tabular}
	\end{small}
	\label{t2}
\end{table*}
	
Table \ref{t1} shows that  existing variance reduced ZO proximal algorithms \cite{huang2019faster,ji2019improved} inevitably use  coordinated ZO  estimator (which is defined in Section \ref{zo_grad}) to estimate full ZO gradient. This is because coordinated ZO  estimator is a more accurate estimation of the true gradient   compared to random ZO estimator, and the corresponding convergence analysis  only needs to handle the error introduced by the smoothing parameter. However, calculating one ZO gradient with coordinated ZO  estimators involves $\mathcal{O}\left( d \right )$ function queries, which makes the FQC of   \cite{huang2019faster,ji2019improved} to be a multiplication between $d$ and $n^{\alpha}$  with $\alpha \geq \frac{1}{2}$. While calculating one full ZO gradient with random ZO estimators involves only $\mathcal{O}\left( 1 \right )$ function queries, which makes it possible to derive variance reduced ZO proximal algorithms with the FQC of summation between  $d$ and $n$. However, it is more challenging  to handle one more error  in addition to the error introduced by the smoothing parameter for non-convex and convex problems.
	
To address this challenging problem, we first proposes a ZO objective decrease (ZOOD) property, which can incorporate not only the error introduced by the smoothing parameter but also the error induced by random ZO estimators in the  upper bound of  convergence rate. Next, we propose two stagewise  reduction frameworks for ZO optimization which automatically derives the convergence results for convex and non-convex problems if the convergence rate of the inner solver satisfies  the ZOOD property. Then, we propose ZOR-ProxSVRG and ZOR-ProxSAGA, two variance reduced ZO proximal  algorithms with fully random ZO estimators and prove that they satisfy the ZOOD property under the strongly convex setting. With the application of two stagewise reduction frameworks on ZOR-ProxSVRG and ZOR-ProxSAGA, we improve the state-of-the-art FQC from $\mathcal{O}\left(\min\{\frac{dn^{1/2}}{\epsilon^2}, \frac{d}{\epsilon^3}\}\right)$ to  $\tilde{\mathcal{O}}\left(\frac{n+d}{\epsilon^2}\right)$ under $d > n^{\frac{1}{2}}$ for non-convex problems, and from $\mathcal{O}\left(\frac{d}{\epsilon^2}\right)$ to $\tilde{\mathcal{O}}\left(n\log\frac{1}{\epsilon}+\frac{d}{\epsilon}\right)$ for convex problems. Experiments on black-box adversarial attack and logistic regression demonstrate that the FQC of our algorithms are significantly lower than  existing state-of-the-art algorithms. To guarantee our reduction algorithms to converge to the same system accuracy of \cite{ghadimi2016mini,huang2019faster,ji2019improved}, we switch to PROX-ZO-SPIDER if the improvement of objective function is smaller than a threshold. We refer the readers to Section \ref{experiment} for more details.

\noindent \textbf{Novelties.} The  main novelties of this paper are summarized as follows:
\begin{itemize}[leftmargin=0.2in]
		
\item Existing reduction frameworks  \cite{allen2016optimal,lin2015universal} require the inner solver to  output $x'$ with input $x_0$ such that $F(x')-\min_{\mathbf{x}}F(\mathbf{x}) \leq \mathcal{K}\left[F(x_0)-\min_{\mathbf{x}}F(\mathbf{x})\right]$, where $0 < \mathcal{K} < 1$. This requirement does not hold for our ZO algorithms due to the two errors introduced by the smoothing parameter and random ZO estimators. To overcome this obstacle, we propose the ZOOD property to incorporate two different errors. Then we conduct the convergence results for our reduction frameworks  based on the ZOOD property.

\item Existing  variance reduced ZO proximal algorithms use coordinated ZO  estimator to calculate full ZO gradient, which makes their query complexities having a multiplication between $d$ and $n^{\alpha}$  with $\alpha \geq \frac{1}{2}$. However, as described above, a better query complexity with a summation  between $d$ and $n$ can be achieved by using random ZO estimator to approximate full ZO gradient. As a result, we propose our two lightweight variance reduced ZO proximal algorithms. Furthermore, by applying our two  reduction frameworks on them we obtain lower FQC for convex and non-convex problems. 
\end{itemize}

\section{Related Work}
\label{related}
\textbf{Black-Box Reduction Techniques.} \cite{allen2016optimal,lin2015universal} proposed two universal stagewise reduction frameworks on first-order algorithms for convex  problems. In each stage, a quadratic reduction term is added to the original optimization problem to ensure the strong convexity holds. \cite{chen2018universal} proposed a universal stagewise reduction framework on first-order algorithms for non-convex optimization problems, where a quadratic reduction term is added similar to \cite{lin2015universal}. All three frameworks work in a black-box manner, \emph{i.e.}, they can be applied to algorithms without knowing algorithms' implementations, but these frameworks are applicable only to first-order algorithms. Our study  falls into developing black-box reduction frameworks on ZO algorithms both for   convex and  non-convex problems.
	
\noindent\textbf{Zeroth-Order Optimization.} Zeroth-order optimization is a classical technique in the optimization community. For smooth problems, \cite{nesterov2017random} proposed zeroth-order gradient descent (ZO-GD) algorithm. Then \cite{ghadimi2013stochastic} proposed its stochastic counterpart ZO-SGD. \cite{lian2016comprehensive} proposed an asynchronous zeroth-order stochastic gradient (ASZO) algorithm for parallel optimization. \cite{gu2018faster} further improved the convergence rate of ASZO by combining variance reduction technique with coordinated ZO  estimator. \cite{liu2018zeroth} proposed ZO-SVRG based algorithms using three different gradient estimators. \cite{fang2018spider} proposed a SPIDER based zeroth-order method named SPIDER-SZO. \cite{ji2019improved} further improved ZO SVRG based and SPIDER based algorithms. \cite{chen2019zo} proposed zeroth-order adaptive momentum method (ZO-AdaMM). \cite{dvurechensky2018accelerated} proposed an accelerated zeroth-order stochastic algorithm and \cite{chen2020accelerated} proposed ZO-Varag which leverages acceleration and variance-reduced technique. For non-smooth regularized problems, \cite{ghadimi2016mini} proposed randomized stochastic projected gradient free (RSPGF) algorithm. \cite{huang2019faster} proposed two zeroth-order variance reduced algorithms with proximal update, ZO-ProxSVRG and ZO-SAGA, and \cite{ji2019improved} proposed a SPIDER based algorithm, PROX-ZO-SPIDER. Note that, the full ZO gradients in these proximal ZO algorithms are estimated in a coordinated manner which would have a high FQC.

\section{Preliminaries}
\label{pre}
For simplicity, we use $||\cdot||$ to denote the Euclidean norm $||\cdot||_2$, and  $\mathbf{x}^*$ to denote the optimal solution for our problem under strong convexity, \emph{i.e.}, $F(\mathbf{x}^*) = \min_{\mathbf{x} \in \mathbb{R}^d}F(\mathbf{x})$. We give some basic definitions as follows. We use the notation $\tilde{\mathcal{O}}(\cdot)$ to hide logarithmic dependence of $d$ and $n$.
	
\begin{definition}
\label{def1}
For function $g : \mathbb{R}^d \rightarrow \mathbb{R}$, we have 
\begin{itemize}[leftmargin=0.2in]
\vspace{-5pt} \item $g$ is $L$-smooth with respect to the Euclidean norm if $g$ has continuous gradients and $\forall \, \mathbf{x}, \mathbf{y} \in \mathbb{R}^d$, it satisfies $|g(\mathbf{y})-g(\mathbf{x})-\langle \nabla g(\mathbf{x}), \mathbf{y} - \mathbf{x}\rangle| \leq \frac{L}{2}||\mathbf{y}-\mathbf{x}||^2$.
			
\vspace{-5pt} \item $g$ is convex if $\forall \, \mathbf{x}, \mathbf{y} \in \mathbb{R}^d$, it satisfies $g(\mathbf{y}) \geq g(\mathbf{x}) + \langle  \nabla g(\mathbf{x}), \mathbf{y} - \mathbf{x}\rangle $.
			
\vspace{-5pt} \item $g$ is $\gamma$-strongly convex if $g(\mathbf{x})-\frac{\gamma}{2}||\mathbf{x}||^2$ is convex, \emph{i.e.}, it satisfies $g(\mathbf{y}) \geq g(\mathbf{x}) + \langle  \nabla g(\mathbf{x}), \mathbf{y} - \mathbf{x}\rangle + \frac{\gamma}{2}||\mathbf{y} - \mathbf{x}||^2$, $\forall \, \mathbf{x}, \mathbf{y} \in \mathbb{R}^d$.
			
\vspace{-5pt} \item $g$ is $\sigma$-weakly convex if $g(\mathbf{x})+\frac{\sigma}{2}||\mathbf{x}||^2$ is convex, \emph{i.e.}, it satisfies $g(\mathbf{y}) \geq g(\mathbf{x}) + \langle  \nabla g(\mathbf{x}), \mathbf{y} - \mathbf{x}\rangle - \frac{\sigma}{2}||\mathbf{y} - \mathbf{x}||^2$, $\forall \, \mathbf{x}, \mathbf{y} \in \mathbb{R}^d$.
		\end{itemize}
	\end{definition}
	
From Definition \ref{def1}, we know that if $g(\mathbf{x})$ is $L$-smooth, then it is $L$-weakly convex. In consequence, for a function $g(\mathbf{x})$ which is both $L$-smooth and $\sigma$-weakly convex, we can assume that $\sigma \leq L$.
	
\begin{definition}[$\epsilon$-Stationary Point]     \label{stationary}
	For problem \eqref{problem}, we call $\mathbf{x} \in \mathbb{R}^d$ a stationary point if $\mathbf{0} \in \left(\nabla f(\mathbf{x}) + \partial r(\mathbf{x})\right)$, where $\partial r(\mathbf{x}) = \{\mathbf{v} \in \mathbb{R}^d | r(\mathbf{y}) \geq r(\mathbf{x}) + \langle \mathbf{v}, \mathbf{y}-\mathbf{x} \rangle, \forall \mathbf{y} \in \mathbb{R}^d\}$ is the subdifferential of $r$ at $\mathbf{x} \in \mathbb{R}^d$. We call $\mathbf{x} \in R^d$ an $\epsilon$-stationary point if $dist\left(\mathbf{0}, \nabla f(\mathbf{x}) + \partial r(\mathbf{x})\right) \leq \epsilon$, where $dist\left(\mathbf{x}, A\right) = \min_{\mathbf{a} \in A} \|\mathbf{x}-\mathbf{a}\|$.
\end{definition}

\begin{definition}[Moreau Envelope and Proximal Mapping]
For any function $F$ and $\lambda > 0$, the following function is called a Moreau envelope of $F$
\begin{equation}
	F_\lambda(\mathbf{x}) = \min_\mathbf{z}\left\{F(\mathbf{z})+\frac{1}{2\lambda}||\mathbf{z}-\mathbf{x}||^2\right\}.
\end{equation}
Further, the optimal solution to the above problem is denoted as
\begin{align*}
    \textrm{Prox}_{\lambda F}(\mathbf{x}) = \arg\min_{\mathbf{z}} \left\{F(\mathbf{z}) + \frac{1}{2\lambda}||\mathbf{z} - \mathbf{x}||^2\right\}
\end{align*}
\end{definition}

\begin{remark}\label{rmk1}
From \cite{chen2018universal}, we know that $F_\lambda(\mathbf{x})$ is smooth for $\lambda < \sigma^{-1}$ if $F$ is $\sigma$-weakly convex. The  gradient of $F_\lambda(\mathbf{x})$ is given by $\nabla F_\lambda(\mathbf{x}) = \frac{\mathbf{x}-\textrm{Prox}_{\lambda F}(\mathbf{x})}{\lambda}$. Further, \cite{chen2018universal} shows that any point $\mathbf{x} \in \mathbb{R}^d$  satisfying $||\nabla F_\lambda(\mathbf{x})|| \leq \epsilon$ is close to a point  that is  $\epsilon$-stationary regarding $F(\mathbf{z})$ in distance of $\mathcal{O}(\epsilon)$.
\end{remark}	
	
\subsection{Assumptions}
Throughout this paper, we will need the following assumptions.
\begin{assum_abbr}
\vspace{-6pt}
\label{a61}
Each $f_i, i = 1, ..., n$ is $L$-smooth.
\end{assum_abbr}

\begin{assum_abbr}
\vspace{-6pt}
\label{a3}
$f$ is $\gamma$-strongly convex.
\end{assum_abbr}

\begin{assum_abbr}
\vspace{-6pt}
\label{a4}
$f$ is convex.
\end{assum_abbr}

\begin{assum_abbr}
\vspace{-6pt}
\label{a1}
$f$ is $\sigma$-weakly convex.
\end{assum_abbr}

\begin{assum_abbr}
\vspace{-6pt}
\label{a2}
$r$ is convex but possibly non-smooth, and the proximal operator of $r$ can be computed efficiently.
\end{assum_abbr}

\subsection{ZO Gradient Estimation}
\label{zo_grad}
\begin{figure}[htbp]
\center
\includegraphics[width=1.0\linewidth]{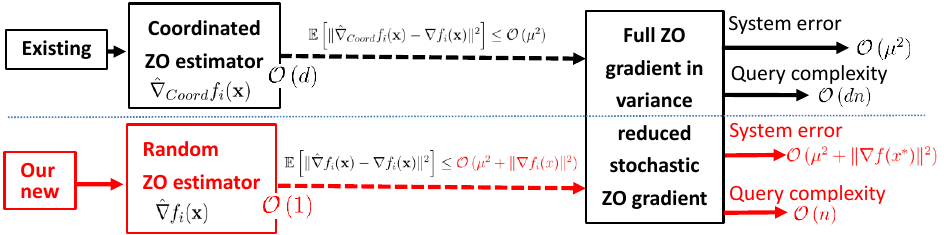}
\vspace*{-7pt}
\caption{Comparison of coordinated and random ZO estimators for variance reduced ZO proximal gradient algorithms.}
\label{fig.sp}
\end{figure}

Given an individual problem $f_i(\mathbf{x})$, we use the two-point  gradient estimator with random direction (abbreviated as \textit{random ZO estimator}) \cite{nesterov2017random,flaxman2004online} to approximate $\nabla f_i(\mathbf{x})$  as follows.
\begin{equation}	\label{FD}
	\hat{\nabla}f_i(\mathbf{x}) = \frac{d}{\mu}[f_i(\mathbf{x}+\mu \mathbf{u}) - f_i(\mathbf{x})]\mathbf{u},
\end{equation}
where $\mu$ is a smoothing parameter, $\mathbf{u} \in \mathbb{R}^d$ is an i.i.d. random direction drawn from a uniform distribution over a unit sphere.	
Note that,  another  commonly used  ZO gradient estimate \cite{liu2018zeroth,ji2019improved,chen2020accelerated}  is \textit{coordinated ZO estimator} of the form of $\label{eqn_coord} \hat{\nabla}_{Coord}f_i(\mathbf{x}) = \sum_{j=1}^d \frac{f_i(\mathbf{x}+\mu \mathbf{e}_j) - f_i(\mathbf{x}-\mu \mathbf{e}_j)}{2\mu}\mathbf{e}_j$, where $\mathbf{e}_j$ is a standard basis vector with 1 at its $j$-th coordinate and $0$ otherwise. As shown in the left part of  Fig. \ref{fig.sp}, the coordinated ZO  estimator is more accurate  than random ZO estimator for approximate ${\nabla}f_i(\mathbf{x})$, while it requires $\mathcal{O}\left( d \right )$ function queries. Conversely, random ZO estimator only requires $\mathcal{O}\left( 1 \right )$ function queries.

\noindent \textbf{Full ZO  Gradient.} No matter which type of variance reduced stochastic ZO gradients  \cite{huang2019faster,ji2019improved} is used,  a full ZO gradient would be involved. For example, \cite{huang2019faster,ji2019improved}  used    \textit{coordinated ZO  estimator} to calculate the full ZO  gradient of the form of $
\frac{1}{n} \sum_{i=1}^n \hat{\nabla}_{Coord}f_i(\mathbf{x})$. As shown in the right part of Fig. \ref{fig.sp}, \textit{coordinated ZO  estimator} would have a smaller system error for variance reduced ZO proximal algorithms, but it has a higher query complexity. In this paper, we will explore using \textit{random ZO estimator} to calculate the full ZO  gradient  of the form of $\hat{\nabla}f(\mathbf{x}) = \frac{1}{n} \sum_{i=1}^n \hat{\nabla}f_i(\mathbf{x})$, and providing the corresponding convergence results to show that it can achieve a lower query complexity.

\section{Zeroth-Order Reduction Frameworks}
\label{ZOOD}
\begin{figure}[htbp]
\center
\includegraphics[width=1.0\linewidth]{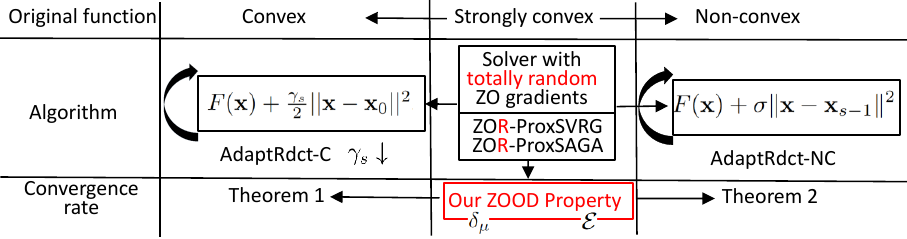}
 \vspace*{-6pt}
 \caption{Principle of our AdaptRdct-C and AdaptRdct-NC.}
 \label{fig.sp2}
 \vspace*{-12pt}
\end{figure}

Due to the benefits of black-box reduction frameworks as discussed previously, we will use  them  to conduct our ZO algorithms to solve \eqref{problem}. However, as discussed before, existing black-box reduction frameworks are limited to first-order algorithms, cannot handle multiple errors induced by random ZO estimators as shown in Fig. \ref{fig.sp}.
To tackle this issue, we first propose a   ZOOD property in Definition \ref{ZOODdef} for our black-box reduction frameworks. Specifically, different to HOOD property proposed by \cite{allen2016optimal} which  requires the objective value to decrease strictly in expectation, that is, $\mathbb{E}\left[F(\mathbf{x}')-\min_{\mathbf{x}}F(\mathbf{x})\right] \leq \mathcal{K}\left[F(\mathbf{x}_0)-\min_{\mathbf{x}}F(\mathbf{x})\right]$, our ZOOD property  absorbs two error terms $\delta_\mu$ and $\mathcal{E}$ which would be corresponding to the error introduced by the smoothing parameter $\mu$ and error introduced by random ZO estimator in our full ZO  gradient, respectively. Thus, our ZOOD property has a relaxed requirement and we will show in Section \ref{apply} that our ZOOD property can be satisfied by our proposed  ZOR-ProxSVRG and ZOR-ProxSAGA with fully random ZO estimators. Based on the ZOOD property, we will propose our two ZO reduction frameworks for convex and non-convex problems respectively. The above principle is also illustrated in Fig. \ref{fig.sp2}.

\begin{definition}[ZOOD Property]
\label{ZOODdef}
We say an algorithm $\mathcal{A}(F, \mathbf{x})$ solving problem \eqref{problem} satisfies the Zeroth-Order Objective Decrease (ZOOD) property with complexity $\mathcal{C}(L, \gamma)$ if, for any starting point $\mathbf{x}_0$, it produces an output $\mathbf{x}' \leftarrow \mathcal{A}(F, \mathbf{x}_0)$ with FQC $\mathcal{C}(L, \gamma)$, such that 
\begin{equation}
	\begin{aligned}
		\mathbb{E}\left[F(\mathbf{x}')-\min_{\mathbf{x}}F(\mathbf{x})-\delta_\mu\right] \leq \mathcal{K}\left[F(\mathbf{x}_0)-\min_{\mathbf{x}}F(\mathbf{x})-\delta_\mu\right] + \mathcal{E},
	\end{aligned}
\end{equation}
where $\delta_\mu$ is a fixed error introduced by the smoothing parameter $\mu$, $0 < \mathcal{K} < 1$ and $\mathcal{E}$ is a constant error.
\end{definition}

\subsection{AdaptRdct-C}
	 \vspace{-3pt}
In this subsection, we  propose a black-box reduction framework  based on the ZOOD property (i.e., \emph{AdaptRdct-C}) for the problem \eqref{problem} in the convex setting.
\emph{AdaptRdct-C} consists of $S$ stages. For each stage,
 we add a quadratic reduction term $\frac{\gamma_s}{2}||\mathbf{x}-\mathbf{x}_0||^2 $ to the original  problem \eqref{problem} such that the problem satisfies the strongly convex property so that we could utilize the linear convergence rate for strongly convex problem to accelerate the original  problem \eqref{problem}. Thus, at the $s$-th stage of \emph{AdaptRdct-C}), our objective function is  the $\gamma_s$-strongly convex problem \eqref{rdctProb1} instead of the original problem \eqref{problem}:
\begin{equation}
\label{rdctProb1}
\begin{aligned}
	F^{(\gamma_s)}(\mathbf{x}) & \overset{\textrm{def}}{=} 
	F(\mathbf{x}) + \frac{\gamma_s}{2}||\mathbf{x}-\mathbf{x}_0||^2 \\ 
	& = \left(f(\mathbf{x}) + \frac{\gamma_s}{2}||\mathbf{x}-\mathbf{x}_0||^2\right) + r(\mathbf{x}).
\end{aligned}
\end{equation} 
Then we run a ZOOD algorithm $\mathcal{A}$ on $F^{(\gamma_s)}(\mathbf{x})$, where the  starting solution  is set as $\mathbf{x}_{s}$, i.e., the output of $\mathcal{A}$ at the $s-1$ stage.   We summarize our \emph{AdaptRdct-C} in Algorithm \ref{algo2}. Note that  the  coefficient  $\gamma_s$ diminishes at an exponential rate. We provide  Theorem \ref{Theorem_3} to show that \emph{AdaptRdct-C} has a linear convergence rate for solving the problem \eqref{problem} in the convex setting. The detailed proof can be found in the supplementary material.

\begin{algorithm}[htb]
\caption{AdaptRdct-C Framework}
	\renewcommand{\algorithmicrequire}{\textbf{Input:}}
	\renewcommand{\algorithmicensure}{\textbf{Output:}}
	\begin{algorithmic}[1]
	\label{algo2}
		\REQUIRE Initial solution  $\mathbf{x}_0$, initial regularization  coefficient $\gamma_0$, epoch number $S$, discount coefficient  $\mathcal{K} \in (0, 1)$, and an inner solver $\mathcal{A}$ satisfying ZOOD  property.
		\FOR {$s = 0, 1, ..., S-1$}
		\STATE $F^{(\gamma_s)}(\mathbf{x}) \overset{\textrm{def}}{=} F(\mathbf{x}) + \frac{\gamma_s}{2}||\mathbf{x}-\mathbf{x}_0||^2$.
		\STATE $\mathbf{x}_{s+1} \leftarrow \mathcal{A}(F^{(\gamma_{s})}, \mathbf{x}_{s})$.
		\STATE $\gamma_{s+1} = \sqrt{\mathcal{K}}\gamma_{s}$.
		\ENDFOR
		\ENSURE $\mathbf{x}_S$.
	\end{algorithmic}
\end{algorithm}

\begin{algorithm}[htb]
\caption{AdaptRdct-NC Framework}
	\renewcommand{\algorithmicrequire}{\textbf{Input:}}
	\renewcommand{\algorithmicensure}{\textbf{Output:}}
	\begin{algorithmic}[1]
	\label{algo3}
		\REQUIRE Initial solution  $\mathbf{x}_0$, epoch number $S$,  $\lambda = \frac{1}{2\sigma}$, and an inner solver $\mathcal{A}$ satisfying ZOOD  property.
		\FOR {$s = 1, 2, ..., S$}
		\STATE $
		\begin{aligned}
			F^{(s)} (\mathbf{x}) & \overset{\textrm{def}}{=} F(\mathbf{x}) + \frac{1}{2\lambda} \|\mathbf{x}-\mathbf{x}_{s-1}\|^2 \\ 
			& = F(\mathbf{x}) + \sigma \|\mathbf{x}-\mathbf{x}_{s-1}\|^2 \\
		\end{aligned}$
		\STATE $\mathbf{x}_s \leftarrow \mathcal{A}(F^{(s)}, \mathbf{x}_{s-1})$.
		\ENDFOR
		\ENSURE $\mathbf{x}_\alpha$ randomly chosen from $\{\mathbf{x}_s\}_{s=1}^{S}$.
	\end{algorithmic}
\end{algorithm}

\begin{theorem}		\label{Theorem_3}
Suppose Assumptions \ref{a61}, \ref{a4} and \ref{a2} are satisfied. Let $\mathbf{x}_0$ be an initial point such that $F(\mathbf{x}_0) - F(\mathbf{x}^*) \leq \Delta$, and $||\mathbf{x}_0 - \mathbf{x}^*||^2 \leq \Theta$. For Algorithm\ref{algo2}, if the inner algorithm $\mathcal{A}$ satisfies the ZOOD property, we have
\begin{equation}
\begin{aligned}
	\mathbb{E}[F(\mathbf{x}_S)-F(\mathbf{x}^*)] 
	\leq  \delta_\mu + \mathcal{K}^S[\Delta - \delta_\mu] + \left (\frac{1}{2}+\frac{2}{\sqrt{\mathcal{K}}} \right )\mathcal{K}^{\frac{S}{2}}\gamma_0\Theta + \frac{\mathcal{E}}{1-\mathcal{K}}.
\end{aligned}
\end{equation}
\end{theorem}

\begin{remark}
The coefficient of the quadratic term in \eqref{rdctProb1}, \emph{i.e.}, $\gamma_s$, diminishes at an exponential rate which contributes to the linear convergence of our AdaptRdct-C. $\delta_\mu$ is the system error induced by the smoothing parameter $\mu$, and $\frac{\mathcal{E}}{1-\mathcal{K}}$ is another   one induced by random ZO estimator.
\end{remark}

\begin{corollary}	\label{coro1}
If the inner algorithm $\mathcal{A}$ satisfies the ZOOD property, the total FQC of \emph{AdaptRdct-C} is $\sum_{s = 0}^{S-1}\mathcal{C}(L, \gamma_s)$ with $S = \mathcal{O}(\log \frac{1}{\epsilon})$, where $\epsilon=\mathbb{E}[F(\mathbf{x}_S)-F(\mathbf{x}^*)] -\delta_\mu -\frac{\mathcal{E}}{1-\mathcal{K}}$.
	\end{corollary}

\subsection{AdaptRdct-NC}
In this subsection, we propose a black-box reduction framework  based on the ZOOD property (i.e., \emph{AdaptRdct-NC}) for solving  problem \eqref{problem} in the non-convex setting.  \emph{AdaptRdct-NC} also consists of multiple stages.  For the $s$-th stage, to make the optimization problem to be strongly convex similar to \emph{AdaptRdct-C}, we add a quadratic term $\|\mathbf{x}-\mathbf{x}_{s-1}\|^2$ with a  coefficient greater than $\sigma/2$ because  the problem \eqref{problem} is  $\sigma$-weakly convex. Different to the exponential decay for the coefficient in \emph{AdaptRdct-C}, the coefficient of the quadratic term here is a constant. Specially, we set the coefficient to be $\sigma$ such that the new problem \eqref{problem} is $\sigma$-strongly convex:
\begin{equation}
\label{rdctProb2}
\begin{aligned}
    F^{(s)}(\mathbf{x}) \overset{\textrm{def}}{=} F(\mathbf{x}) + \sigma ||\mathbf{x}-\mathbf{x}_{s-1}||^2 
    = \left(f(\mathbf{x}) + \sigma ||\mathbf{x}-\mathbf{x}_{s-1}||^2\right) + r(\mathbf{x}).
\end{aligned}
\end{equation}
Then we call a zeroth-order algorithm $\mathcal{A}$ which satisfies the ZOOD property to solve $F^{(s)}(\mathbf{x})$, where the  starting point is $\mathbf{x}_{s-1}$ (i.e., the output of $\mathcal{A}$ at the $s-1$ stage). We summarize our \emph{AdaptRdct-NC} in Algorithm \ref{algo3}. The convergence property of \emph{AdaptRdct-NC} is given in Theorem \ref{Theorem_4}, where the proof can be found in the supplementary material.

\begin{theorem}		\label{Theorem_4}
Suppose Assumptions  \ref{a61}, \ref{a1} and \ref{a2} are satisfied. Let $\mathbf{x}_0$ be an initial solution  such that $F(\mathbf{x}_0) - F(\mathbf{x}^*) \leq \Delta$.  For Algorithm \ref{algo3}, if the inner algorithm $\mathcal{A}$ satisfies the ZOOD property, we have
\begin{align}
\label{t3eq8}
    \mathbb{E}[||\nabla F_\lambda(\mathbf{x}_\alpha)||^2] \leq \frac{8\sigma(2\mathcal{K}+1)\Delta}{(1-\mathcal{K})(S+1)} + 24\sigma\delta_\mu + \frac{24\mathcal{E}}{1-\mathcal{K}}.
\end{align}
\end{theorem}
\begin{remark}
As mentioned in Remark \ref{rmk1}, a point $\mathbf{x} \in \mathbb{R}^d$  with $||\nabla F_\lambda(\mathbf{x})|| \leq \epsilon$ is actually an $\epsilon$-stationary point of $F$. Thus, the convergence analysis of \emph{AdaptRdct-NC} is built on  $\|\nabla f_\lambda (\mathbf{x})\|^2$ as shown in Theorem \ref{Theorem_4}. In the right side of inequality of \eqref{t3eq8}, $24\sigma\delta_\mu$ and $\frac{24\mathcal{E}}{1-\mathcal{K}}$ correspond to the  system errors induced by the smoothing parameter $\mu$ and  random ZO estimator respectively.
\end{remark}
	
\begin{corollary}	\label{coro2}
    If the inner algorithm $\mathcal{A}$ satisfies the ZOOD property,  the total FQC of \emph{AdaptRdct-NC} is $S \times \mathcal{C}(L, \sigma)$ with $S = \mathcal{O}(\frac{1}{\epsilon^2})$, where $\epsilon^2=\mathbb{E}||\nabla F(\mathbf{x}_\alpha)||^2 - 24\sigma\delta_\mu -\frac{24\mathcal{E}}{1-\mathcal{K}}$.
\end{corollary}

\section{Applications on Lightweight  Variance Reduced ZO Proximal  Algorithms}
\label{apply}
In this section, we propose our two lightweight variance reduced ZO proximal algorithms (i.e.,  ZOR-ProxSVRG and ZOR-ProxSAGA) by fully utilizing random ZO estimators. Then we  provide the theoretical results (i.e., Theorems \ref{remark1} and \ref{remark2}) to show that the ZOOD property which is the key requirement of our \emph{AdaptRdct-C} and \emph{AdaptRdct-NC} hold for our  ZOR-ProxSVRG and ZOR-ProxSAGA under the strongly convex setting. Finally, we  apply \emph{AdaptRdct-C} and \emph{AdaptRdct-NC} on  ZOR-ProxSVRG and ZOR-ProxSAGA and obtain the corresponding FQC according to Theorems \ref{Theorem_3} and \ref{Theorem_4}. The theoretical results show that our FQC on convex problems outperforms the existing one (Corollary \ref{corollary6}), and  our FQC on non-convex problems outperforms the ones of vanilla ZO-ProxSVRG, ZO-ProxSAGA and state-of-the-art PROX-ZO-SPIDER  when $d > n^{\frac{1}{6}}$ or  $d > n^{\frac{1}{2}}$ (Corollary \ref{corollary7}).

\begin{algorithm}[htb]
\caption{ZOR-ProxSVRG Algorithm}
\renewcommand{\algorithmicrequire}{\textbf{Input:}}
\renewcommand{\algorithmicensure}{\textbf{Output:}}
\begin{algorithmic}[1]
	\label{algo1}
	\REQUIRE Stepsize $\eta$, update frequency $m$, batch size $b$, epoch $S$, starting point $\mathbf{x}_m^0 = \tilde{\mathbf{x}}_0 \in \mathbb{R}^d$.
	\FOR {$s = 1, 2, ..., S$}
	\STATE 
	$\mathbf{x}_0^s = \tilde{\mathbf{x}}_{s-1}$;
	\STATE Compute $\hat{\mathbf{v}}_s = \hat{\nabla}f(\tilde{\mathbf{x}}_{s-1}) $.
	\FOR {$k = 0, 1, ..., m-1$}
	\STATE Uniformly randomly choose $\mathcal{I}_k \subseteq \{1, 2, ..., n\}$ with $|\mathcal{I}_k| = b$ and $\mathbf{u}_{i_k}^s$ from a unit sphere for $i_k \in \mathcal{I}_k$.
	\STATE $\hat{\mathbf{v}}_k^s = \hat{\nabla}f_{\mathcal{I}_k}(\mathbf{x}_k^s) - \hat{\nabla}f_{\mathcal{I}_k}(\tilde{\mathbf{x}}_{s-1}) + \hat{\mathbf{v}}_s $.
	\STATE $\mathbf{x}_{k+1}^s = \textrm{Prox}_{\eta r}\left(\mathbf{x}_k^s - \eta \hat{\mathbf{v}}_k^s\right)$.
	\ENDFOR
	\STATE Set $\tilde{\mathbf{x}}_s = \mathbf{x}_k^s$ for randomly chosen $ k \in \{0, 1, ..., m-1\}$.
	\ENDFOR
	\ENSURE $\mathbf{x}_m^S$.
\end{algorithmic}
\end{algorithm}

\begin{algorithm}[htb]
\caption{ZOR-ProxSAGA Algorithm}
	\renewcommand{\algorithmicrequire}{\textbf{Input:}}
	\renewcommand{\algorithmicensure}{\textbf{Output:}}
\begin{algorithmic}[1]
\label{algo4}
	\REQUIRE Stepsize $\eta$, batch size $b$, epoch $K$, starting point $\mathbf{x}^0 \in R^d$, auxiliary vectors $\{\phi_i^0\}_{i=1}^n = \mathbf{x}^0$ for each $i$.
	\STATE Compute $\hat{\mathbf{g}}^0 = \hat{\nabla} f(\mathbf{x}^0) $.
	\FOR {$k = 0, 1, 2, ..., K$}
	\STATE Uniformly randomly choose $\mathcal{I}_k \subseteq \{1, 2, ..., n\}$ with $|\mathcal{I}_k| = b$ and $\mathbf{u}_{i_k}^s$ from a unit sphere for $i_k \in \mathcal{I}_k$.
	\STATE $\hat{\mathbf{v}}^k = \hat{\nabla}f_{\mathcal{I}_k}(\mathbf{x}^k) - \hat{\nabla}f_{\mathcal{I}_k}(\mathbf{\phi}_{i}^k) + \hat{\mathbf{g}}^k $.
	\STATE $\mathbf{x}^{k+1} = \textrm{Prox}_{\eta r}\left(\mathbf{x}^k - \eta \hat{\mathbf{v}}^k\right)$.
	\STATE $\mathbf{\phi}_{i_k}^{k+1} = \mathbf{x}^k$.
	\STATE $\hat{\mathbf{g}}^{k+1} = \hat{\mathbf{g}}^k + \frac{b}{n} (\hat{\nabla}f_{\mathcal{I}_k}(\mathbf{x}^k) - \hat{\nabla}f_{\mathcal{I}_k}(\mathbf{\phi}_{i}^k))$.
	\ENDFOR
	\ENSURE $\mathbf{x}^{K}$.
\end{algorithmic}
\end{algorithm}

\subsection{Lightweight Variance Reduced ZO Proximal Algorithms}
In this part, we  use random  ZO estimator to lighten variance reduced ZO proximal algorithms and propose our ZOR-ProxSVRG and ZOR-ProxSAGA as follows.
	  
\noindent \textbf{ZOR-ProxSVRG:}  
Xiao and Zhang proposed ProxSVRG in \cite{xiao2014proximal}. The key step of ProxSVRG is to generate  a snapshot (denoted as $\tilde{\mathbf{x}}$) of $\mathbf{x}$ after a certain number  of iterations, and the full gradient at $\tilde{\mathbf{x}}$ is used to build a modified stochastic gradient estimation, which is a gradient blending 	 $\mathbf{v} = \nabla f_i(\mathbf{x}) - \nabla f_i(\tilde{\mathbf{x}}) + \nabla f(\tilde{\mathbf{x}}) $, where $\mathbf{v}$ denotes the gradient estimate at $\mathbf{x}$, $i \in [n]$ is chosen uniformly randomly, and $\nabla f(\tilde{\mathbf{x}}) = \frac{1}{n}\sum_{j=1}^n \nabla f_j(\tilde{\mathbf{x}})$. In  ZO setting, the true gradient is replaced by gradient estimator and we sample a mini batch of $i$ instead of a single one. The ZOR-ProxSVRG algorithm is described in Algorithm \ref{algo1}. Note that all ZO gradients are estimated by a \textit{random ZO estimator}. We provide Theorem \ref{remark1} to show that  our ZOR-ProxSVRG satisfies the ZOOD property.

\begin{theorem}
\label{remark1}
	Suppose Assumptions \ref{a61}, \ref{a3} and \ref{a2} hold,  ZOR-ProxSVRG satisfies the ZOOD property with  $\mathcal{C}(L, \gamma) = \mathcal{O}\left(n + \frac{dL}{\gamma}\right)$, $\mathcal{K} = \frac{\beta_2}{\beta_1}$, $\mathcal{E} = \mathcal{O}\left(\|\nabla f(\mathbf{x}^*)\|^2\right)$ and $\delta_\mu = \frac{2\eta L\mu^2}{\beta_1-\beta_2}\left(1+(\frac{3}{b}+1)\eta Ld^2\right)$, where $\beta_1 = 2\eta\left(1-\frac{24dL\eta}{b}\right)$ and $\beta_2 = \frac{2}{m\gamma}+\frac{48dL\eta^2}{bm} + 16dL\eta^2\left(\frac{3}{b}+2\right)$.
\end{theorem}
	
\noindent \textbf{ZOR-ProxSAGA:} \ 
Defazio et al. proposed first-order SAGA algorithm in \cite{defazio2014saga}. The key step of SAGA is to keep a table of gradients of previous results $\nabla f(\mathbf{\phi}_{i\in [n]})$. The modified stochastic gradient estimate is $\mathbf{v} = \nabla f_i(\mathbf{x}) - \nabla f_i(\mathbf{\phi}_i) + g(\mathbf{\phi})$, where $\mathbf{v}$ denotes the gradient estimate at $\mathbf{x}$, $i \in [n]$ is chosen uniformly randomly, and $g(\mathbf{\phi}) = \frac{1}{n}\sum_{j=1}^n \nabla f_j(\mathbf{\phi}_j)$. Under ZO setting, the true gradient is replaced by gradient estimator (\textit{i.e.}, \textit{random ZO estimator } in this paper) and we sample a mini batch of $i$ instead of a single one. The ZOR-ProxSAGA algorithm is described in Algorithm \ref{algo4}. We provide Theorem \ref{remark2} to show that  our ZOR-ProxSAGA satisfies the ZOOD property.

\begin{theorem}
\label{remark2}
Suppose Assumptions  \ref{a61}, \ref{a3} and \ref{a2} hold, ZOR-ProxSAGA satisfies the ZOOD property with $\mathcal{C}(L, \gamma) = \mathcal{O}\left(n + \frac{dL}{\gamma}\log\left(\max\left\{dL, n\gamma\right\}\right)\right)$, $\mathcal{E} = \mathcal{O}\left(d\|\nabla f(\mathbf{x}^*)\|^2\right)$ and $\delta_\mu = \frac{2(b+3)L^2d^2\mu^2}{b\gamma}$. 
\end{theorem}

\subsection{Applying Reduction Frameworks on ZOR-ProxSVRG and ZOR-ProxSAGA}
With the above results,  we can present the FQC of  \emph{AdaptRdct-C} and \emph{AdaptRdct-NC} variants of  ZOR-ProxSVRG and ZOR-ProxSAGA in Corollaries \ref{corollary6} and \ref{corollary7} respectively.

Before providing the theoretical results, we first define $G_\emph{C}^2$ and $G_\emph{NC}^2$ in Definition \ref{a6} which are the upper bounds of $\|\nabla f^{(\gamma_s)}(\mathbf{x}_s^*)||^2$ and $\|\nabla f^{(s)}(\mathbf{x}_s^*)\|^2$, respectively,  involved in our reduction frameworks.	
\vspace{-5pt}
\begin{definition}[$G_\emph{C}^2$ and $G_\emph{NC}^2$]
\label{a6}
	Define $G_\emph{C}^2$ to be the number such that $\|\nabla f^{(\gamma_s)}(\mathbf{x}_s^*)||^2 \leq G_\emph{C}^2\;$ for each stage $s = 0, 1, ..., S-1$ of Algorithm \ref{algo2}, where $f^{(\gamma_s)}(\mathbf{x}) = f(\mathbf{x}) + \frac{\gamma_s}{2}\|\mathbf{x}-\mathbf{x}_0\|^2$ and $\mathbf{x}_s^* = \arg\min_\mathbf{x} F^{(\gamma_s)}(\mathbf{x})$.  $G_\emph{NC}^2$ is  the number such that $\|\nabla f^{(s)}(\mathbf{x}_s^*)\|^2 \leq G_\emph{NC}^2\;$   for each stage $s = 1, 2, ..., S$ of Algorithm \ref{algo3}, where $f^{(s)}(\mathbf{x}) = f(\mathbf{x}) + \sigma\|\mathbf{x}-\mathbf{x}_{s-1}\|^2$ and $\mathbf{x}_s^* = \arg\min_\mathbf{x} F^{(s)}(\mathbf{x})$.
\end{definition}
	
\begin{corollary}
\label{corollary6}
Suppose Assumptions \ref{a61}, \ref{a4} and \ref{a2}  are satisfied, applying \emph{AdaptRdct-C} on ZOR-ProxSVRG and ZOR-ProxSAGA. From Corollary \ref{coro1}, Theorem \ref{remark1} and \ref{remark2}, denote $G_*^2 = \max\left\{G_\emph{C}^2, G_\emph{NC}^2\right\}$, we have that the total FQC to solve problem \eqref{problem} is $\tilde{\mathcal{O}}\left(n\log\frac{1}{\epsilon}+\frac{d}{\epsilon}\right)$ with  $\epsilon=\mathbb{E}[F(\mathbf{x}_S)-F(\mathbf{x}^*)] -\mathcal{O}\left( \mu^2 + G_*^2\right)$. 
\end{corollary}

\begin{remark}
From Table \ref{t2} and Corollary \ref{corollary6}, we know that \emph{AdaptRdct-C} improves the state-of-the-art FQC for convex problem from   $\mathcal{O}\left(\frac{d}{\epsilon^2}\right)$ to $\tilde{\mathcal{O}}\left(n\log\frac{1}{\epsilon}+\frac{d}{\epsilon}\right)$  before reaching the  accuracy of $\mathcal{O}\left( \mu^2 + G_*^2\right)$.
\end{remark}

\begin{corollary}
\label{corollary7}
Suppose Assumptions \ref{a61}, \ref{a1} and \ref{a2}  are satisfied, applying \emph{AdaptRdct-NC} on ZOR-ProxSVRG and ZOR-ProxSAGA. From Corollary \ref{coro2}, Theorem \ref{remark1} and \ref{remark2}, denote $G_*^2 = \max\left\{G_\emph{C}^2, G_\emph{NC}^2\right\}$, we have that the total FQC for  solving problem \eqref{problem} is $\tilde{\mathcal{O}}\left(\frac{n+d}{\epsilon^2}\right)$ with $\epsilon^2=\mathbb{E}||\nabla F(\mathbf{x}_\alpha)||^2 - \mathcal{O}\left( \mu^2 + G_*^2\right)$.
\end{corollary}

\begin{remark}
    Under the non-convex setting, from Table \ref{t1} and Corollary \ref{corollary7}, we know that the FQC of \emph{AdaptRdct-NC} outperforms the one of PROX-ZO-SPIDER \cite{ji2019improved}  if $d > n^{\frac{1}{2}}$, and the ones of ZO-ProxSVRG and ZO-ProxSAGA if $d > n^{\frac{1}{6}}$, before achieving the  accuracy of $\mathcal{O}\left( \mu^2 + G_*^2\right)$. This indicates the strengths of our proposed methods in solving high dimensional optimization problems, which could be ubiquitous in modern machine learning and deep learning.
\end{remark}

\begin{figure*}[htb]
\vspace{-20pt}
\centering
\subfigure[Cifar-10]{
\includegraphics[width = 0.31\textwidth, height = 3cm]{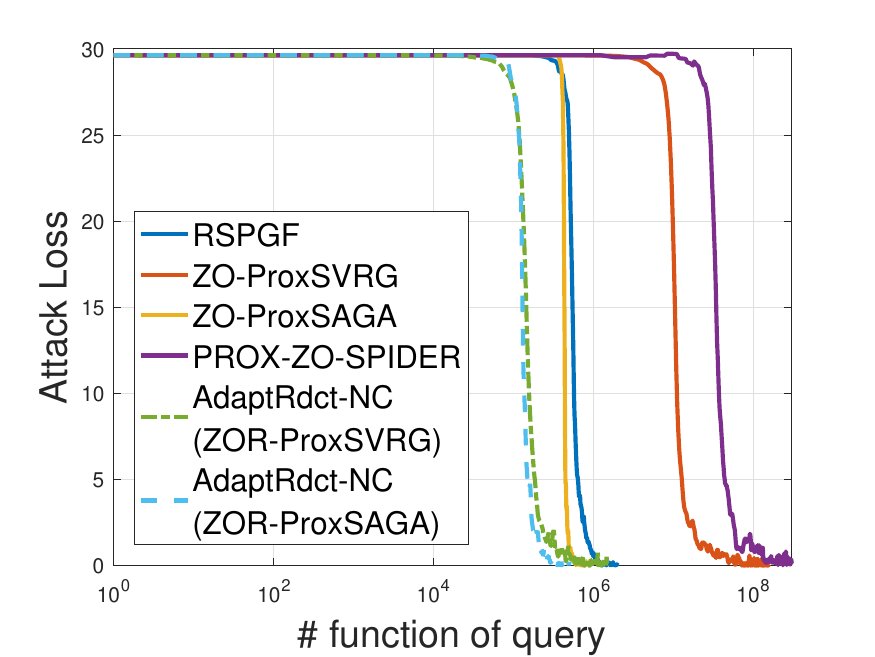}
}
\subfigure[fmnist]{
\includegraphics[width = 0.31\textwidth, height = 3cm]{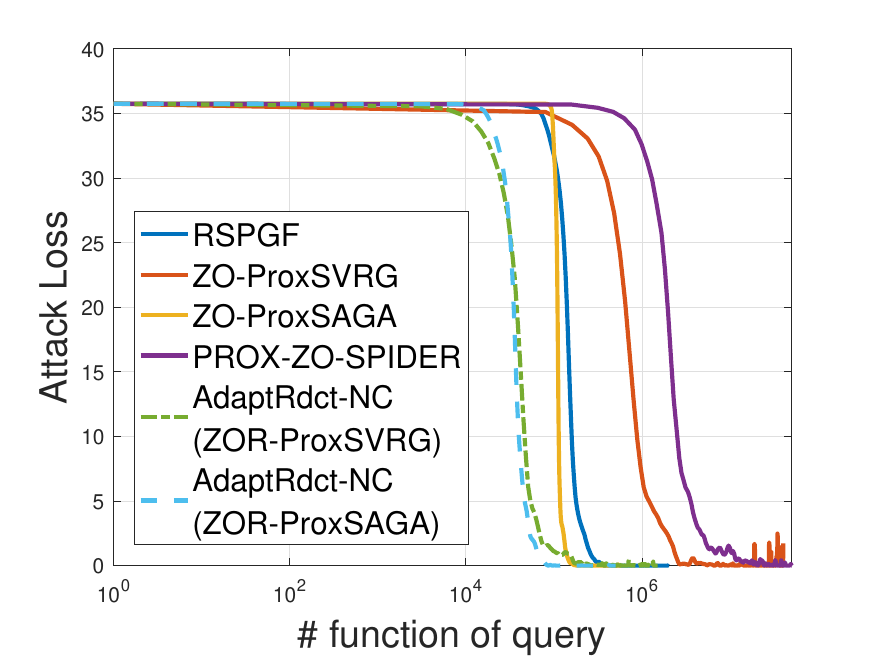}
}
\subfigure[Mnist]{
\includegraphics[width = 0.31\textwidth, height = 3cm]{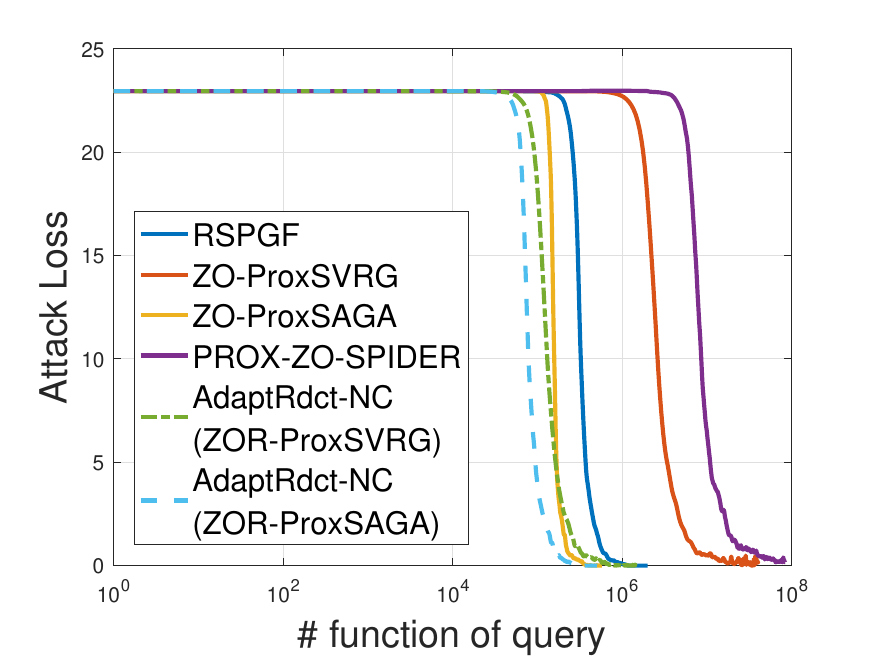}
}
\vspace{-8pt}
\caption{Comparison of black-box attack methods on three well-trained DNNs, with $\sigma = $1e-3.}
\label{adver}
\vspace{-12pt}
\end{figure*}

\begin{figure*}[htb]
\centering
\subfigure[A9a, $\gamma_0 = $5e-4]{
\includegraphics[width = 0.31\textwidth, height = 3cm]{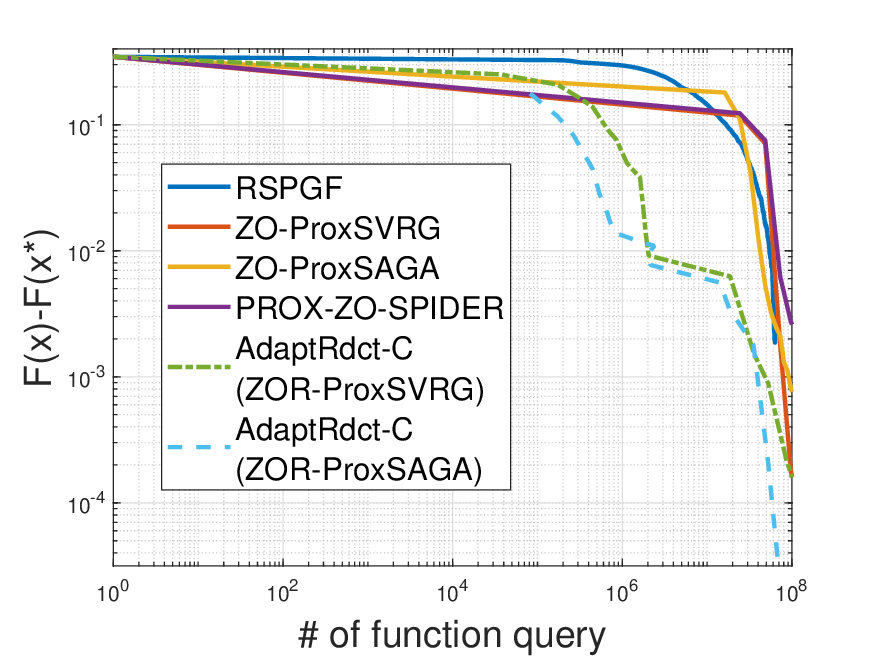}
\label{a9a_convex}
}
\subfigure[Ijcnn1, $\gamma_0 = $5e-4]{
\includegraphics[width = 0.31\textwidth, height = 3cm]{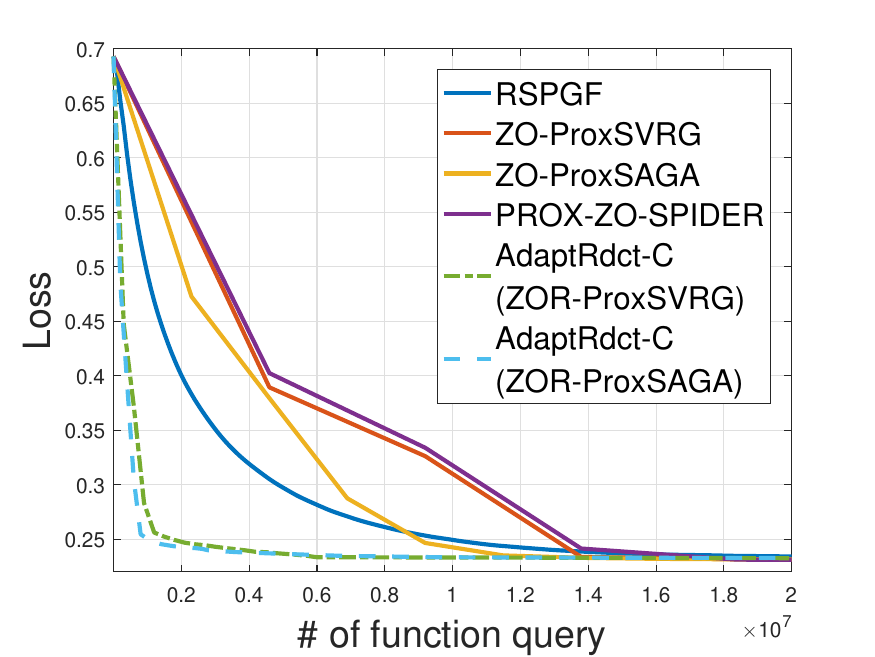}
}
\subfigure[QSAR, $\gamma_0 = $5e-3]{
\includegraphics[width = 0.31\textwidth, height = 3cm]{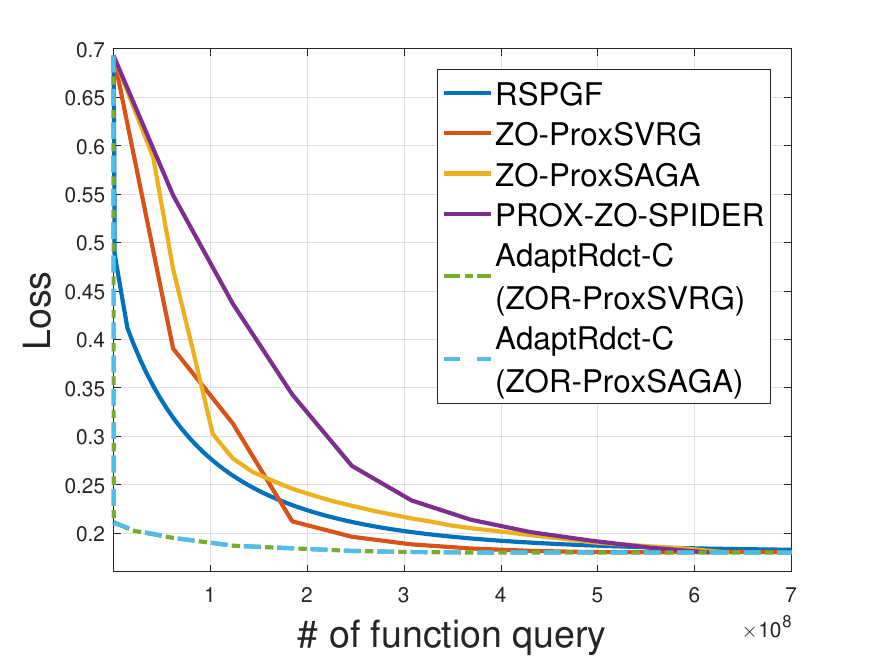}
\label{qsar_convex}
}
\\
\vspace{-12pt}
\subfigure[A9a, $\sigma = $5e-4]{\label{a9a_non_convex}
\includegraphics[width = 0.31\textwidth, height = 3cm]{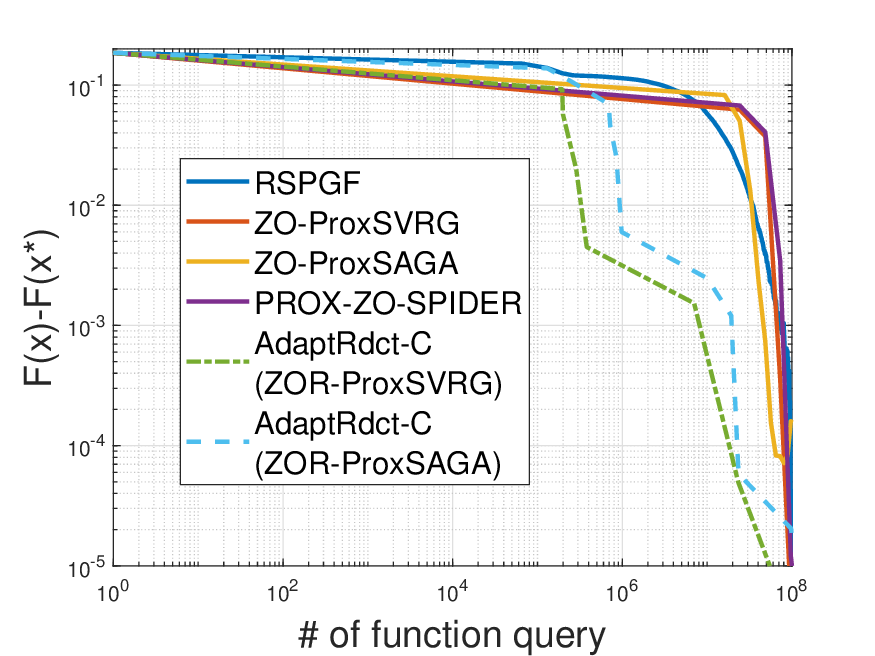}
}
\subfigure[Ijcnn1, $\sigma = $5e-4]{
\includegraphics[width = 0.31\textwidth, height = 3cm]{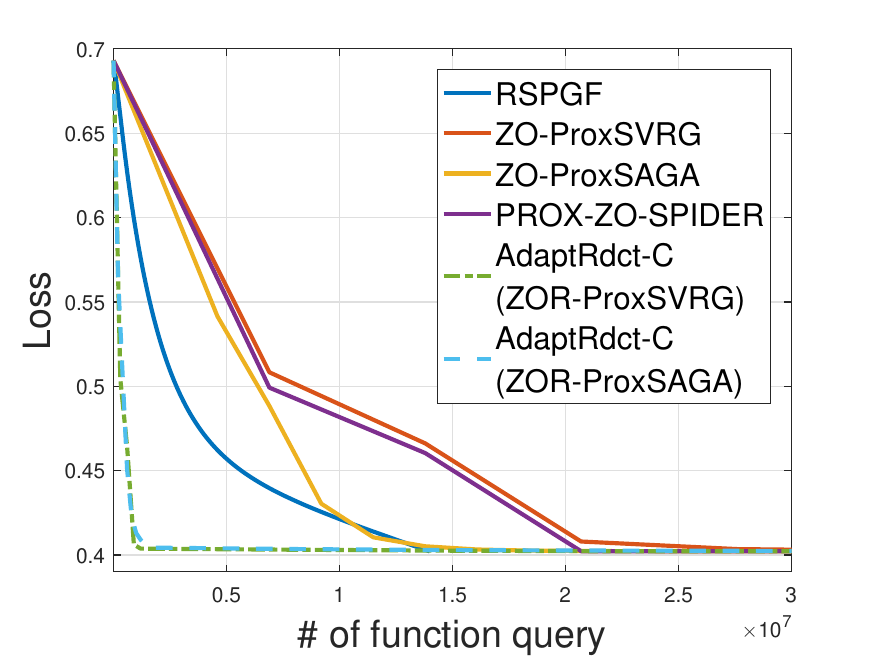}
}
\subfigure[QSAR, $\sigma = $5e-3]{
\includegraphics[width = 0.31\textwidth, height = 3cm]{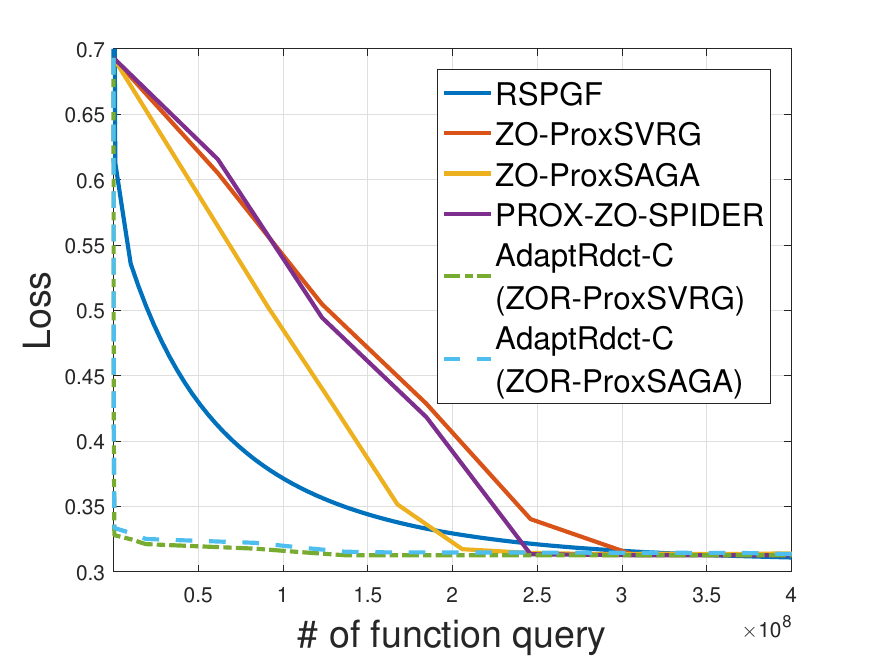}
\label{qsar_non_convex}
}
\\
\vspace{-8pt}
\caption{Comparison of different ZO algorithms for logistic regression problems. (a)-(c) Convex. (d)-(f) Non-convex. (a) and (d) are plotted with residue error (i.e., $F(x)-F(x^*)$) in the log-scale.}
\label{logistic}
\vspace{-12pt}
\end{figure*}

\section{Experimental Results}
\label{experiment}
In this section, we verify the effectiveness of our ZO reduction frameworks with two experiments on several real-world datasets.
The first experiment is generating  black-box adversarial examples which is a non-convex problem. The second experiment addresses logistic regression under convex and non-convex settings, respectively. All algorithms in Table \ref{t1} are the  compared algorithms. Note that our \emph{AdaptRdct-C} and \emph{AdaptRdct-NC} variants  of ZOR-ProxSVRG and ZOR-ProxSAGA may have larger system errors than the compared algorithms as shown in Table \ref{t1}. The error $G_\emph{NC}^2$ is small enough for the first experiment to generate black-box adversarial examples  successfully. For the second experiment, to guarantee our algorithms to converge to the same system accuracy of \cite{ghadimi2016mini,huang2019faster,ji2019improved}, we switch our algorithms to PROX-ZO-SPIDER if the improvement of objective function (\textit{i.e.}, $F(\mathbf{x}_s)-F(\mathbf{x}_{s-1})$) is smaller than a threshold. We set the threshold to 1e-3 and 3e-4 for the convex and the non-convex logistic regression, respectively.

\subsection{Generation of Black-Box Adversarial Examples}
In image classification, adversary attack crafts input images with imperceptive perturbation to mislead a trained classifier. The resulting perturbed images are called adversarial examples, which are commonly used to understand the robustness of learning models. Under the black-box setting, the attackers only have access to the function value. It is obvious that this problem falls into the framework of ZO optimization.
	
For the target black-box model, we choose three well-trained DNNs $\mathcal{F}(\cdot) = \left[\mathcal{F}_1(\cdot), ..., \mathcal{F}_K(\cdot)\right]$, where $\mathcal{F}_k(\cdot)$ returns the prediction score of the $k$-th class. They are trained on three datasets, \emph{i.e.},  Cifar-10, fashion Mnist (fmnist) and Mnist. For each model, we attack a batch of $n$ correctly-classified images $\{\mathbf{a}_i\}_{i=1}^n$ from the same class, and adopt the following black-box attacking loss. The $i$-th individual loss function $f_i(\mathbf{x}) + r(\mathbf{x})$ is given by
\begin{equation}
\max \left( \mathcal{F}_{y_i}\left(\mathbf{a}_i\right) - \mathcal{F}_{y_{tar}}\left(\mathbf{a}_i+\mathbf{x}\right), 0\right) + \lambda_1 \|\mathbf{x}\|_1 + \lambda_2 \|\mathbf{x}\|_2^2,
\end{equation}
where $\mathbf{x}$ is the learned adversarial perturbation, and $\mathbf{x} \in [0, 1]$. In the experiment, we choose $n = 50$, with a batch-size of 10. $\lambda_1$ and $\lambda_2$ are set to 1e-3 and 1 respectively.
The number of directions we use to construct the random gradient estimator \eqref{FD} is 10. We compare our reduction algorithms with four algorithms: RSPGF, ZO-ProxSVRG, ZO-ProxSAGA and PROX-ZO-SPIDER. Our reduction algorithms and RSPGF use the random ZO estimator  and the other all use coordinated ZO  estimator. Fig. \ref{adver} shows the convergence of attack loss, the attacks are successful when the attack loss is close to 0. The results show that ZOR-ProxSVRG and ZOR-ProxSAGA with \emph{AdaptRdct-NC} have a better performance with respect to FQC  than all other algorithms.

\vspace{-5pt}
\subsection{Convex and Nonconvex Logistic Regression}
\vspace{-3pt}
In this subsection, we mainly consider logistic regression and its variant. To conduct experiments on \emph{AdaptRdct-C}, we first choose the classical logistic regression problem $f(\mathbf{x}) = \frac{1}{n}\sum_{i=1}^n \ell(\mathbf{w}^T\mathbf{x}_i, y_i)$ and $F(\mathbf{x}) = f(\mathbf{x}) + \lambda_1 \|\mathbf{x}\|_1 + \lambda_2 \|\mathbf{x}\|_2^2$, where $\mathbf{x}_i \in \mathbb{R}^d$ denote the features, $y_i \in \left\{ 0, 1\right\}$ are the labels and $\ell(\cdot)$ is the cross-entropy loss. It is obvious that $f$ is convex. For non-convex problems, we add a non-convex regularizer $\sum_{i=1}^d \mathbf{w}_i^2/(1 + \mathbf{w}_i^2)$ to the convex problem $f(\mathbf{x})$. Then we get $\bar{f}(\mathbf{x}) = \frac{1}{n}\sum_{i=1}^n \ell(\mathbf{w}^T\mathbf{x}_i, y_i) + \alpha\sum_{i=1}^d \mathbf{w}_i^2/(1 + \mathbf{w}_i^2)$ and $\bar{F}(\mathbf{x}) = \bar{f}(\mathbf{x}) + \lambda_1 \|\mathbf{x}\|_1 + \lambda_2 \|\mathbf{x}\|_2^2$. In the experiment, we set $\lambda_1 = $1e-3 and $\lambda_2 = $1e-5. For these problems, we conduct the experiments on  two LIBSVM datasets   \cite{chang2011libsvm}, \emph{i.e.}, the a9a  ($n = 32,561$, $d = 123$), Ijcnn1 ($n = 49,990$, $d = 22$) datasets, and one UCI dataset, the QSAR dataset\footnote{The link of QSAR dataset is  \url{https://archive.ics.uci.edu/ml/datasets/QSAR+oral+toxicity\#?tdsourcetag=s_pctim_aiomsg}.} ($n = 8,992, d = 1,024$). We compare our reduction algorithms with four algorithms: RSPGF, ZO-ProxSVRG, ZO-ProxSAGA and PROX-ZO-SPIDER. Our reduction algorithms and RSPGF use the random ZO estimator  and the others all use coordinated ZO estimator. Fig. \ref{logistic} shows that the FQCs of ZOR-ProxSVRG and ZOR-ProxSAGA with \emph{AdaptRdct-C} and \emph{AdaptRdct-NC} are significantly lower the ones of all  existing algorithms even for low dimensional data.

\section{Conclusion}

In this paper, we proposed two reduction frameworks to lower the complexity of zeroth-order proximal algorithms for convex and non-convex problems, respectively. Our frameworks work in a black-box manner, thus they can be applied to a wide range of algorithms without knowing their implementations. To the best of our knowledge, our frameworks are the first to lower the FQC of ZO algorithms. Specifically, we studied the performance of applying our frameworks on two variance reduced proximal gradient algorithms using fully random ZO estimators, \textit{i.e.} ZOR-ProxSVRG and ZOR-ProxSAGA. The results indicate that our approach has a significant lower FQC than existing algorithms before achieving the system accuracy. Our theoretic study and experimental results highlight the advantages of combining zeroth-order optimization with the reduction techniques.

\section*{Appendix A: Proof of Theorem  \ref{Theorem_3}}
\begin{custheorem}1		
Suppose Assumptions \ref{a61}, \ref{a4} and \ref{a2} are satisfied. Let $\mathbf{x}_0$ be a starting vector such that $F(\mathbf{x}_0) - F(\mathbf{x}^*) \leq \Delta$, and $||\mathbf{x}_0 - \mathbf{x}^*||^2 \leq \Theta$. For Algorithm \ref{algo2}, we have
\begin{equation}
\mathbb{E}[F(\mathbf{x}_S)-F(\mathbf{x}^*)] \leq \delta_\mu + \mathcal{K}^S[\Delta - \delta_\mu] + \left (\frac{1}{2}+\frac{2}{\sqrt{\mathcal{K}}} \right )\mathcal{K}^{\frac{S}{2}}\gamma_0\Theta + \frac{\mathcal{E}}{1-\mathcal{K}}
\end{equation}
\end{custheorem}

\begin{proof}
Denote $\mathbf{x}_s^* = \arg\min_{\mathbf{x}}F^{(\gamma_s)}(\mathbf{x})$. By the strong convexity of $F^{(\gamma_s)}(\mathbf{x})$, we have
\[
\mathbb{E}\left[F^{(\gamma_s)}(\mathbf{x}_{s}^*)-F^{(\gamma_s)}(\mathbf{x}^*)\right] \leq -\frac{\gamma_s}{2}\mathbb{E}\left[||\mathbf{x}_s^* - \mathbf{x}^*||^2\right]
\]
Using the fact that $F^{(\gamma_s)}(\mathbf{x}_s^*) \geq F(\mathbf{x}_s^*)$, as well as the definition $F^{(\gamma_s)}(\mathbf{x}^*) = F(\mathbf{x}^*) + \frac{\gamma_s}{2}||\mathbf{x}^*-\mathbf{x}_0||^2$, we immediately have
\begin{equation}
\mathbb{E}\left[F(\mathbf{x}_s^*) - F(\mathbf{x}^*) - \frac{\gamma_s}{2}||\mathbf{x}^*-\mathbf{x}_0||^2\right] \leq -\frac{\gamma_s}{2}\mathbb{E}\left[||\mathbf{x}_s^* - \mathbf{x}^*||^2\right]
\end{equation}
Rearranging the terms, we get
\begin{equation}	\label{t3eq1}
\frac{\gamma_s}{2}||\mathbf{x}_0 - \mathbf{x}^*||^2 - \frac{\gamma_s}{2}\mathbb{E}\left[||\mathbf{x}_s^* - \mathbf{x}^*||^2\right] \geq \mathbb{E}\left[F(\mathbf{x}_s^*)-F(\mathbf{x}^*)\right] \geq 0
\end{equation}
Thus we have
\begin{equation}	\label{t3eq2}
\mathbb{E}\left[\|\mathbf{x}_s^*-\mathbf{x}^*\|^2\right] \leq ||\mathbf{x}_0 - \mathbf{x}^*||^2
\end{equation}
	
The ZOOD property of $\mathcal{A}$ ensures that
\begin{equation}	\label{t3eq3}
\mathbb{E}\left[F^{(\gamma_s)}(\mathbf{x}_{s+1})-F^{(\gamma_s)}(\mathbf{x}_s^*)-\delta_\mu\right] \leq \mathcal{K}\mathbb{E}\left[F^{(\gamma_s)}(\mathbf{x}_{s})-F^{(\gamma_s)}(\mathbf{x}_s^*)-\delta_\mu\right] + \mathcal{E}
\end{equation}
	
Define $D_s = \mathbb{E}\left[F^{(\gamma_s)}(\mathbf{x}_s) - F^{(\gamma_s)}(\mathbf{x}_s^*)\right]$. At epoch 0, we have upper bound $D_0 = F^{(\gamma_0)}(\mathbf{x}_0) - F^{(\gamma_0)}(\mathbf{x}_0^*) \leq F(\mathbf{x}_0) - F(\mathbf{x}^*)$. For each epoch $s\geq1$, we compute that
\begin{equation}	\label{t3eq4}
\begin{aligned}
&D_s = \mathbb{E}\left[F^{(\gamma_s)}(\mathbf{x}_s) - F^{(\gamma_s)}(\mathbf{x}_s^*)\right] \\
\overset{\textrm{\ding{172}}}{=}& \mathbb{E}\left[F^{(\gamma_{s-1})}(\mathbf{x}_s) - \frac{\gamma_{s-1}-\gamma_{s}}{2}||\mathbf{x}_s-\mathbf{x}_0||^2\right] - \mathbb{E}\left[F^{(\gamma_{s-1})}(\mathbf{x}_s^*) - \frac{\gamma_{s-1}-\gamma_{s}}{2}||\mathbf{x}_s^*-\mathbf{x}_0||^2\right] \\
\overset{\textrm{\ding{173}}}{\leq}& \mathbb{E}\left[F^{(\gamma_{s-1})}(\mathbf{x}_s) - \frac{\gamma_{s-1}-\gamma_{s}}{2}||\mathbf{x}_s-\mathbf{x}_0||^2\right] \\ 
&-\mathbb{E}\left[F^{(\gamma_{s-1})}(\mathbf{x}_{s-1}^*) + \frac{\gamma_{s-1}}{2}||\mathbf{x}_s^* - \mathbf{x}_{s-1}^*||^2\right] + \frac{\gamma_{s-1}-\gamma_s}{2}\mathbb{E}\left[||\mathbf{x}_s^*-\mathbf{x}_0||^2\right] \\
\leq& \mathbb{E}\left[F^{(\gamma_{s-1})}(\mathbf{x}_s)-F^{(\gamma_{s-1})}(\mathbf{x}_{s-1}^*)\right] + \frac{\gamma_{s-1}-\gamma_s}{2}\mathbb{E}\left[||\mathbf{x}_s^*-\mathbf{x}_0||^2\right] \\
\overset{\textrm{\ding{174}}}{\leq}& \mathbb{E}\left[F^{(\gamma_{s-1})}(\mathbf{x}_s)-F^{(\gamma_{s-1})}(\mathbf{x}_{s-1}^*)\right] + (\gamma_{s-1}-\gamma_{s})\mathbb{E}\left[||\mathbf{x}_s^* - \mathbf{x}^*||^2 + ||\mathbf{x}_0 - \mathbf{x}^*||^2\right] \\
\overset{\textrm{\ding{175}}}{\leq}& \mathbb{E}\left[F^{(\gamma_{s-1})}(\mathbf{x}_s)-F^{(\gamma_{s-1})}(\mathbf{x}_{s-1}^*)\right] + 2(1-\sqrt{\mathcal{K}})\gamma_{s-1}||\mathbf{x}_0-\mathbf{x}^*||^2 \\
\end{aligned}
\end{equation}
Above, \ding{172} uses the definition of $F^{(\gamma_s)}(\mathbf{x})$; \ding{173} uses the convexity of $F^{(\gamma_{s-1})}(\mathbf{x})$ and the fact that $\mathbf{x}_{s-1}^*$ is the minimizer of $F^{(\gamma_{s-1})}(\mathbf{x})$; \ding{174} uses the inequality $||\mathbf{a}-\mathbf{b}||^2 \leq 2||\mathbf{a}||^2 + 2||\mathbf{b}||^2$; \ding{175} follows from (\ref{t3eq2}) and $\gamma_s = \sqrt{\mathcal{K}}\gamma_{s-1}$. That implies
\begin{equation}	\label{t3eq5}
\begin{aligned}
D_s - \delta_\mu &\leq \mathbb{E}[F^{(\gamma_{s-1})}(\mathbf{x}_s)-F^{(\gamma_{s-1})}(\mathbf{x}_{s-1}^*)-\delta_\mu] + 2(1-\sqrt{\mathcal{K}})\gamma_{s-1}||\mathbf{x}_0-\mathbf{x}^*||^2 \\
&\leq \mathcal{K}(D_{s-1}-\delta_\mu) + 2(1-\sqrt{\mathcal{K}})\gamma_{s-1}||\mathbf{x}_0-\mathbf{x}^*||^2 + \mathcal{E} \\
\end{aligned}
\end{equation}
The second inequality follows from (\ref{t3eq3}).
	
Recursively applying the above inequality, we have
\begin{equation}	\label{t3eq6}
\begin{aligned}
&D_S - \delta_\mu \\
\leq& \mathcal{K}^S(D_0 - \delta_\mu) + ||\mathbf{x}_0-\mathbf{x}^*||^2\cdot 2\cdot (1-\sqrt{\mathcal{K}})[\gamma_{S-1} + \mathcal{K}\gamma_{S-2} + \mathcal{K}^2\gamma_{S-3} + ... + \mathcal{K}^{S-1}\gamma_0] + \mathcal{E}\sum_{i=0}^{S-1}\mathcal{K}^i \\
=& \mathcal{K}^S(D_0 - \delta_\mu) + ||\mathbf{x}_0-\mathbf{x}^*||^2\cdot 2\cdot (1-\sqrt{\mathcal{K}})\gamma_{S-1}[1+\sqrt{\mathcal{K}}+\mathcal{K}+ ... + \sqrt{\mathcal{K}}^{S-1}] + \frac{1-\mathcal{K}^S}{1-\mathcal{K}}\mathcal{E} \\
\leq& \mathcal{K}^S(D_0 - \delta_\mu) + ||\mathbf{x}_0-\mathbf{x}^*||^2\cdot 2\cdot \gamma_{S-1} + \frac{\mathcal{E}}{1-\mathcal{K}} \\
=& \mathcal{K}^S(D_0 - \delta_\mu) + ||\mathbf{x}_0-\mathbf{x}^*||^2\frac{2}{\sqrt{\mathcal{K}}}\gamma_{S} + \frac{\mathcal{E}}{1-\mathcal{K}} \\
\end{aligned}
\end{equation}
The first inequality and the last inequality follows from the update rule that $\gamma_{s} = \sqrt{\mathcal{K}}\gamma_{s-1}$
	
Finally, we obtain
\begin{equation}	\label{t3eq7}
\begin{aligned}
&\mathbb{E}[F(\mathbf{x}_S)-F(\mathbf{x}^*) - \delta_\mu] \\
\leq& \mathbb{E}[F^{(\gamma_S)}(\mathbf{x}_S)-F^{(\gamma_S)}(\mathbf{x}^*) - \delta_\mu +\frac{\gamma_S}{2}||\mathbf{x}_0-\mathbf{x}^*||^2] \\
\leq& \mathbb{E}[F^{(\gamma_S)}(\mathbf{x}_S)-F^{(\gamma_S)}(\mathbf{x}_S^*) - \delta_\mu +\frac{\gamma_S}{2}||\mathbf{x}_0-\mathbf{x}^*||^2] \\
\leq& \mathcal{K}^S[F(\mathbf{x}_0)-F(\mathbf{x}^*) - \delta_\mu] + (\frac{1}{2}+\frac{2}{\sqrt{\mathcal{K}}})\gamma_S||\mathbf{x}_0-\mathbf{x}^*||^2 + \frac{\mathcal{E}}{1-\mathcal{K}} \\
\leq& \mathcal{K}^S[F(\mathbf{x}_0)-F(\mathbf{x}^*) - \delta_\mu] + (\frac{1}{2}+\frac{2}{\sqrt{\mathcal{K}}})\mathcal{K}^{\frac{S}{2}}\gamma_0||\mathbf{x}_0-\mathbf{x}^*||^2 + \frac{\mathcal{E}}{1-\mathcal{K}} \\
\end{aligned}
\end{equation}
The first inequality uses the definition of $F^{(\gamma_s)}(\mathbf{x})$ and the fact that $F(\mathbf{x}) \leq F^{(\gamma_s)}(\mathbf{x})$; the second inequality uses the fact that $\mathbf{x}_S^*$ is the minimizer of $F^{(\gamma_s)}(\mathbf{x}_S)$; the third inequality uses (\ref{t3eq6}), and the last inequality uses the update rule $\gamma_{s} = \sqrt{\mathcal{K}}\gamma_{s-1}$. Note that $F(\mathbf{x}_0)-F(\mathbf{x}^*) \leq \Delta$ and $||\mathbf{x}_0-\mathbf{x}^*||^2 \leq \Theta$, then we have
\begin{equation}
\mathbb{E}[F(\mathbf{x}_S)-F(\mathbf{x}^*)] 
\leq \delta_\mu + \mathcal{K}^S[\Delta - \delta_\mu] + (\frac{1}{2}+\frac{2}{\sqrt{\mathcal{K}}})\mathcal{K}^{\frac{S}{2}}\gamma_0\Theta + \frac{\mathcal{E}}{1-\mathcal{K}} \\
\end{equation}
Then we finish the proof.
\end{proof}

\section*{Appendix B: Proof of Theorem  \ref{Theorem_4}}

\begin{custheorem}2		
Suppose Assumptions \ref{a61}, \ref{a1} and \ref{a2} are satisfied. Let $\mathbf{x}_0$ be a starting vector such that $F(\mathbf{x}_0) - F(\mathbf{x}^*) \leq \Delta$. By running Algorithm \ref{algo3}, we have
\begin{align}
\mathbb{E}[||\nabla F_\lambda(\mathbf{x}_\alpha)||^2] \leq \frac{8\sigma(2\mathcal{K}+1)\Delta}{(1-\mathcal{K})(S+1)} + 24\sigma\delta_\mu + \frac{24\mathcal{E}}{1-\mathcal{K}}
\end{align}
\end{custheorem}

\begin{proof}
With the definition of Moreau envelope, we have $F_\lambda(\mathbf{x}_{s-1}) = \min_\mathbf{x}F^{(s)}(\mathbf{x})$. Denote $\mathbf{x}_s^* = \arg\min_\mathbf{x} F^{(s)}(\mathbf{x}) = Prox_{\lambda F}(\mathbf{x}_{s-1})$, we get $\nabla F_\lambda(\mathbf{x}_{s-1}) = \frac{\mathbf{x}_{s-1}-\mathbf{x}_s^*}{\lambda}$.
Then with the definition of $F^{(s)}(\mathbf{x})$, we have
\begin{equation}	\label{t3eq9}
F(\mathbf{x}_{s-1}) = F^{(s)}(\mathbf{x}_{s-1}) \geq F^{(s)}(\mathbf{x}_s^*) = F(\mathbf{x}_s^*) + \frac{1}{2\lambda}||\mathbf{x}_{s-1} - \mathbf{x}_s^*||^2
\end{equation}
The ZOOD property ensures that 
\begin{equation}	\label{t3eq10}
\mathbb{E}\left[F^{(s)}(\mathbf{x}_s) - F^{(s)}(\mathbf{x}_s^*) - \delta_\mu\right] \leq \mathcal{K}\left[F^{(s)}(\mathbf{x}_{s-1}) - F^{(s)}(\mathbf{x}_s^*) - \delta_\mu\right] + \mathcal{E}
\end{equation}
which implies
\begin{equation}
\mathbb{E}\left[F^{(s)}(\mathbf{x}_s)\right] \leq F^{(s)}(\mathbf{x}_s^*) + \mathcal{K}\left[F^{(s)}(\mathbf{x}_{s-1})-F^{(s)}(\mathbf{x}_s^*)\right] - (\mathcal{K}-1)\delta_\mu + \mathcal{E} \\
\end{equation}
Thus we have
\begin{equation}	\label{t3eq11}
\begin{aligned}
&\mathbb{E}\left[F(\mathbf{x}_s) + \frac{1}{2\lambda}||\mathbf{x}_s-\mathbf{x}_{s-1}||^2\right] = \mathbb{E}\left[F^{(s)}(\mathbf{x}_s)\right] \\
\leq& F^{(s)}(\mathbf{x}_s^*) + \mathcal{K}\left[F^{(s)}(\mathbf{x}_{s-1})-F^{(s)}(\mathbf{x}_s^*)\right] - (\mathcal{K}-1)\delta_\mu + \mathcal{E} \\
\leq& F(\mathbf{x}_{s-1}) + \mathcal{K}\left[F^{(s)}(\mathbf{x}_{s-1})-F^{(s)}(\mathbf{x}_s^*)\right] - (\mathcal{K}-1)\delta_\mu + \mathcal{E} \\
\end{aligned}
\end{equation}
The last inequality follows from (\ref{t3eq9}). On the other hand, we have
\begin{equation}	\label{t3eq12}
\begin{aligned}
||\mathbf{x}_s - \mathbf{x}_{s-1}||^2 &= ||\mathbf{x}_s - \mathbf{x}_s^* + \mathbf{x}_s^* - \mathbf{x}_{s-1}||^2 \\
&= ||\mathbf{x}_s - \mathbf{x}_s^*||^2 + ||\mathbf{x}_s^* - \mathbf{x}_{s-1}||^2 + 2\left\langle\mathbf{x}_s - \mathbf{x}_s^*, \mathbf{x}_s^* - \mathbf{x}_{s-1}\right\rangle \\
&\geq -||\mathbf{x}_s - \mathbf{x}_s^*||^2 + \frac{1}{2}||\mathbf{x}_{s-1} - \mathbf{x}_s^*||^2
\end{aligned}
\end{equation}
The inequality follows from Young's inequality. Combining with (\ref{t3eq11}), we then get
\begin{equation}	\label{t3eq13}
\begin{aligned}
&\mathbb{E}\left[\frac{1}{4\lambda}||\mathbf{x}_{s-1} - \mathbf{x}_s^*||^2\right] \leq \mathbb{E}\left[\frac{1}{2\lambda}||\mathbf{x}_s - \mathbf{x}_{s-1}||^2 + \frac{1}{2\lambda}||\mathbf{x}_s - \mathbf{x}_s^*||^2\right] \\
\leq& \mathbb{E}\left[F(\mathbf{x}_{s-1})-F(\mathbf{x}_s) + \mathcal{K}[F^{(s)}(\mathbf{x}_{s-1})-F^{(s)}(\mathbf{x}_s^*)] - (\mathcal{K}-1)\delta_\mu + \mathcal{E} + \frac{1}{2\lambda}||\mathbf{x}_s - \mathbf{x}_s^*||^2\right] \\
\overset{\textrm{\ding{172}}}{\leq}& \mathbb{E}\left[F(\mathbf{x}_{s-1})-F(\mathbf{x}_s) + \mathcal{K}[F^{(s)}(\mathbf{x}_{s-1})-F^{(s)}(\mathbf{x}_s^*)] - (\mathcal{K}-1)\delta_\mu + \mathcal{E}\right] \\
&+  \mathbb{E}\left[\frac{1}{\lambda(\lambda^{-1}-\sigma)}[F^{(s)}(\mathbf{x}_s)-F^{(s)}(\mathbf{x}_s^*)]\right] \\
\overset{\textrm{\ding{173}}}{\leq}& \mathbb{E}\left[F(\mathbf{x}_{s-1})-F(\mathbf{x}_s) + \frac{2-\lambda\sigma}{1-\lambda\sigma}\left(\mathcal{K}\left[F^{(s)}(\mathbf{x}_{s-1})-F^{(s)}(\mathbf{x}_s^*)\right] - (\mathcal{K}-1)\delta_\mu + \mathcal{E}\right)\right] \\
\overset{\textrm{\ding{174}}}{=}& \mathbb{E}\left[F(\mathbf{x}_{s-1})-F(\mathbf{x}_s) + 3\left(\mathcal{K}\left[F^{(s)}(\mathbf{x}_{s-1})-F^{(s)}(\mathbf{x}_s^*)\right] - (\mathcal{K}-1)\delta_\mu + \mathcal{E}\right)\right] \\
\end{aligned}
\end{equation}
\ding{172} holds because $F^{(s)}(\mathbf{x})$ is $(\lambda^{-1}-\sigma)$-strongly convex, \ding{173} holds due to (\ref{t3eq10}), and \ding{174} holds due to $\lambda = \frac{1}{2\sigma}$. Next, we bound $F^{(s)}(\mathbf{x}_{s-1})-F^{(s)}(\mathbf{x}_s^*)$ given that $\mathbf{x}_{s-1}$ is fixed. According to the definition of $F^{(s)}(\mathbf{x})$, we have
\begin{equation}	\label{t3eq14}
\begin{aligned}
F^{(s)}(\mathbf{x}_{s-1})-F^{(s)}(\mathbf{x}_s^*) &= F^{(s)}(\mathbf{x}_s)-F^{(s)}(\mathbf{x}_s^*) + F^{(s)}(\mathbf{x}_{s-1})-F^{(s)}(\mathbf{x}_s) \\
&= F^{(s)}(\mathbf{x}_s)-F^{(s)}(\mathbf{x}_s^*) + F(\mathbf{x}_{s-1})-F(\mathbf{x}_s)-\frac{1}{2\lambda}||\mathbf{x}_s-\mathbf{x}_{s-1}||^2 \\
&\leq F^{(s)}(\mathbf{x}_s)-F^{(s)}(\mathbf{x}_s^*) + F(\mathbf{x}_{s-1})-F(\mathbf{x}_s)
\end{aligned}
\end{equation}
Taking expectation over randomness in the $s$-th stage on both sides, we have
\begin{equation}	\label{t3eq15}
\begin{aligned}
F^{(s)}(\mathbf{x}_{s-1})-F^{(s)}(\mathbf{x}_s^*) &\leq \mathbb{E}\left[F^{(s)}(\mathbf{x}_s)-F^{(s)}(\mathbf{x}_s^*)\right] + \mathbb{E}\left[F(\mathbf{x}_{s-1})-F(\mathbf{x}_s)\right] \\
&\leq \mathcal{K}\left[F^{(s)}(\mathbf{x}_{s-1}) - F^{(s)}(\mathbf{x}_s^*)\right] - (\mathcal{K}-1)\delta_\mu + \mathcal{E} + \mathbb{E}[F(\mathbf{x}_{s-1})-F(\mathbf{x}_s)]
\end{aligned}
\end{equation}
The second inequality follows from (\ref{t3eq10}). 
Rearranging the terms, we have
\begin{equation}	\label{t3eq16}
F^{(s)}(\mathbf{x}_{s-1})-F^{(s)}(\mathbf{x}_s^*) \leq \frac{1}{1-\mathcal{K}} \mathbb{E}[F(\mathbf{x}_{s-1})-F(\mathbf{x}_s)] + \delta_\mu + \frac{\mathcal{E}}{1-\mathcal{K}}
\end{equation}
Plug this upper bound into (\ref{t3eq13}), then we have
\begin{equation}	\label{t3eq17}
\begin{aligned}
\mathbb{E}\left[\frac{1}{4\lambda}||\mathbf{x}_{s-1} - \mathbf{x}_s^*||^2\right] &\leq \mathbb{E}\left[F(\mathbf{x}_{s-1})-F(\mathbf{x}_s) + \frac{3\mathcal{K}}{1-\mathcal{K}}[F(\mathbf{x}_{s-1})-F(\mathbf{x}_s)]+3\delta_\mu + \frac{3\mathcal{E}}{1-\mathcal{K}}\right] \\
&= \frac{2\mathcal{K}+1}{1-\mathcal{K}}\mathbb{E}\left[F(\mathbf{x}_{s-1})-F(\mathbf{x}_s)\right] + 3\delta_\mu + \frac{3\mathcal{E}}{1-\mathcal{K}} \\
\end{aligned}
\end{equation}
Since $||\nabla F_\lambda(\mathbf{x}_{s-1})|| = \frac{||\mathbf{x}_{s-1} - \mathbf{x}_s^*||}{\lambda}$, we have
\begin{equation}	\label{t3eq18}
\frac{\lambda}{4}\mathbb{E}\left[||\nabla F_\lambda(\mathbf{x}_{s-1})||^2\right] \leq \frac{2\mathcal{K}+1}{1-\mathcal{K}}\mathbb{E}\left[F(\mathbf{x}_{s-1})-F(\mathbf{x}_s)\right] + 3\delta_\mu + \frac{3\mathcal{E}}{1-\mathcal{K}}
\end{equation}
By summing over $s = 1, 2, ..., S+1$, we get
\begin{equation}	\label{t3eq19}
\frac{\lambda}{4}\sum_{s=1}^{S+1}\mathbb{E}\left[||\nabla F_\lambda(\mathbf{x}_{s-1})||^2\right] \leq \frac{2\mathcal{K}+1}{1-\mathcal{K}}\mathbb{E}\left[F(\mathbf{x}_0)-F(\mathbf{x}_{S+1})\right] + 3(S+1)\delta_\mu + \frac{3(S+1)\mathcal{E}}{1-\mathcal{K}}
\end{equation}
With the choice of $\mathbf{x}_\alpha$, we have
\begin{equation}	\label{t3eq20}
\frac{\lambda}{4}(S+1)\mathbb{E}[||\nabla F_\lambda(\mathbf{x}_\alpha)||^2] \leq \frac{2\mathcal{K}+1}{1-\mathcal{K}}[F(\mathbf{x}_0)-F(\mathbf{x}_{S+1})] + 3(S+1)\delta_\mu + \frac{3(S+1)\mathcal{E}}{1-\mathcal{K}}
\end{equation}
Since $\lambda = \frac{1}{2\sigma}$, we have
\begin{equation}	\label{t3eq21}
\mathbb{E}[||\nabla F_\lambda(\mathbf{x}_\alpha)||^2] \leq \frac{8\sigma(2\mathcal{K}+1)}{(1-\mathcal{K})(S+1)}[F(\mathbf{x}_0)-F(\mathbf{x}_{S+1})] + 24\sigma\delta_\mu + \frac{24\mathcal{E}}{1-\mathcal{K}}
\end{equation}
Note that $F(\mathbf{x}_0)-F(\mathbf{x}_{S+1}) \leq F(\mathbf{x}_0)-F(\mathbf{x}^*) = \Delta$. 
Thus, we get Theorem \ref{Theorem_4}.
\end{proof}

\section*{Appendix C: Proof of Theorem \ref{remark1}}

Firstly we present one additional lemma (i.e. Lemma \ref{lemma5}) for ZOR-ProxSVRG, then we present  Theorem \ref{Theorem_1} for our ZOR-ProxSVRG.

We introduce a smoothing function. Define 
\begin{equation}
\label{f_mu}
f_\mu (\mathbf{x}) = \mathbb{E}_{\mathbf{u}\sim \mathbf{U}_b} [f(\mathbf{x}+\mu\mathbf{u})],
\end{equation}
where $\mathbf{U}_b$ is a uniform distribution over the unit Euclidean ball.

\begin{cuslemma}1	\label{lemma5}
Suppose Assumptions  \ref{a61}, \ref{a4} and \ref{a2} hold. We have
\begin{equation}
\begin{aligned}
&\mathbb{E}\left[\|\hat{\mathbf{v}}_k^s-\nabla f_\mu(\mathbf{x}_k^s)\|^2\right] \\
\leq& \frac{24dL}{b} \left[F(\mathbf{x}_k^s)-F(\mathbf{x}^*)\right] + 8dL\left(\frac{3}{b}+2\right)\left[F(\tilde{\mathbf{x}}_{s-1})-F(\mathbf{x}^*)\right] + 8d\|\nabla f(\mathbf{x}^*)\|^2 + \left(\frac{3}{b}+1\right)L^2d^2\mu^2
\end{aligned}
\end{equation}
\end{cuslemma}
\begin{proof}
Since $\hat{\mathbf{v}}_k^s = \hat{\nabla}f_{\mathcal{I}_k}(\mathbf{x}_k^s) - \hat{\nabla}f_{\mathcal{I}_k}(\tilde{\mathbf{x}}_{s-1}) + \hat{\mathbf{g}}_s$, we have
\begin{equation}
\begin{aligned}
&\mathbb{E}\left[\|\hat{\mathbf{v}}_k^s-\nabla f_\mu(\mathbf{x}_k^s)\|^2\right] \\
=& \mathbb{E}\left[\|\hat{\nabla}f_{\mathcal{I}_k}(\mathbf{x}_k^s) - \hat{\nabla}f_{\mathcal{I}_k}(\tilde{\mathbf{x}}_{s-1}) + \hat{\mathbf{v}}_s-\nabla f_\mu(\mathbf{x}_k^s)\|^2\right] \\
=& \mathbb{E}\left[\|\hat{\nabla}f_{\mathcal{I}_k}(\mathbf{x}_k^s) - \hat{\nabla}f_{\mathcal{I}_k}(\tilde{\mathbf{x}}_{s-1}) + \nabla f_\mu(\tilde{\mathbf{x}}_{s-1}) - \nabla f_\mu(\mathbf{x}_k^s) + \hat{\mathbf{v}}_s-\nabla f_\mu(\tilde{\mathbf{x}}_{s-1})\|^2\right] \\
\overset{\textrm{\ding{172}}}{\leq}& 2\mathbb{E}\left[\|\hat{\nabla}f_{\mathcal{I}_k}(\mathbf{x}_k^s) - \hat{\nabla}f_{\mathcal{I}_k}(\tilde{\mathbf{x}}_{s-1})-\mathbb{E}\left[\hat{\nabla}f_{\mathcal{I}_k}(\mathbf{x}_k^s) - \hat{\nabla}f_{\mathcal{I}_k}(\tilde{\mathbf{x}}_{s-1})\right]\|^2\right] \\
&+ 2\mathbb{E}\left[\|\hat{\nabla}f(\tilde{\mathbf{x}}_{s-1}) - \nabla f_\mu(\tilde{\mathbf{x}}_{s-1})\|^2\right] \\
\overset{\textrm{\ding{173}}}{\leq}&  2\mathbb{E}\left[\|\hat{\nabla}f_{\mathcal{I}_k}(\mathbf{x}_k^s) - \hat{\nabla}f_{\mathcal{I}_k}(\tilde{\mathbf{x}}_{s-1})\|^2\right] + 4d\|\nabla f(\tilde{\mathbf{x}}_{s-1})\|^2 + L^2d^2\mu^2 \\
\overset{\textrm{\ding{174}}}{\leq}& \frac{2}{bn}\sum_{i = 1}^{n} \mathbb{E}\left[\|\hat{\nabla}f_{i}(\mathbf{x}_k^s; \mathbf{u}_k^s) - \hat{\nabla}f_{i}(\tilde{\mathbf{x}}_{s-1}; \mathbf{u}_k^s)\|^2\right] + 4d\|\nabla f(\tilde{\mathbf{x}}_{s-1})-\nabla f(\mathbf{x}^*) + \nabla f(\mathbf{x}^*)\|^2 + L^2d^2\mu^2 \\
\overset{\textrm{\ding{175}}}{\leq}& \frac{12d}{bn}\sum_{i = 1}^{n}\mathbb{E}\left[\|\nabla f_i(\mathbf{x}_k^s)- \nabla f_i(\mathbf{x}^*)\|^2+\|\nabla f_i(\tilde{\mathbf{x}}_{s-1})- \nabla f_i(\mathbf{x}^*)\|^2\right] + \frac{3L^2d^2\mu^2}{b} \\
&+ 8d\|\nabla f(\tilde{\mathbf{x}}_{s-1})-\nabla f(\mathbf{x}^*)\|^2 + 8d\|\nabla f(\mathbf{x}^*)\|^2 + L^2d^2\mu^2\\
\overset{\textrm{\ding{176}}}{\leq}& \frac{24dL}{b} \left[F(\mathbf{x}_k^s)-F(\mathbf{x}^*)+F(\tilde{\mathbf{x}}_{s-1})-F(\mathbf{x}^*)\right] + \frac{3L^2d^2\mu^2}{b} \\
&+ \frac{8d}{n}\sum_{i=1}^n\|\nabla f_i(\tilde{\mathbf{x}}_{s-1})-\nabla f_i(\mathbf{x}^*)\|^2 + 8d\|\nabla f(\mathbf{x}^*)\|^2 + L^2d^2\mu^2\\
\overset{\textrm{\ding{177}}}{\leq}& \frac{24dL}{b} \left[F(\mathbf{x}_k^s)-F(\mathbf{x}^*)\right] + 8dL\left(\frac{3}{b}+2\right)\left[F(\tilde{\mathbf{x}}_{s-1})-F(\mathbf{x}^*)\right] + 8d\|\nabla f(\mathbf{x}^*)\|^2 + \left(\frac{3}{b}+1\right)L^2d^2\mu^2,
\end{aligned}
\end{equation}
where \ding{172} comes from Young's inequality; \ding{173} comes from the fact that $\mathbb{E}\left[\|\mathbf{a}-\mathbb{E}\left[\mathbf{a}\right]\|^2\right] \leq \mathbb{E}\left[\|\mathbf{a}\|^2\right]$ and Lemma \ref{lemma1}; \ding{174} directly comes from Lemma \ref{lemma6}; \ding{175} comes from Lemma \ref{lemma3} with $\beta = 1$ and Young's inequality; \ding{176} comes from Lemma \ref{lemma4} and Jensen's inequality; \ding{177} comes from Lemma \ref{lemma4} and rearranging the terms. Then we complete the proof.
\end{proof}

\begin{custheorem}5	\label{Theorem_1}
Suppose Assumptions \ref{a61}, \ref{a3} and \ref{a2} hold. By running Algorithm \ref{algo1}, we have
\begin{equation}
\mathbb{E}\left[F(\tilde{\mathbf{x}}_S)-F(\mathbf{x}^*)\right] \leq \delta_\mu + \left(\frac{\beta_2}{\beta_1}\right)^S\left[F(\tilde{\mathbf{x}}_{0})-F(\mathbf{x}^*) - \delta_\mu\right] + \frac{16\eta^2d\|\nabla f(\mathbf{x}^*)\|^2}{\beta_1-\beta_2},
\end{equation}
where $\beta_1 = 2\eta\left(1-\frac{24dL\eta}{b}\right), \beta_2 = \frac{2}{m\gamma}+\frac{48dL\eta^2}{bm} + 16dL\eta^2\left(\frac{3}{b}+2\right)$, $\delta_\mu = \frac{2\eta L\mu^2}{\beta_1-\beta_2}\left(1+(\frac{3}{b}+1)\eta Ld^2\right)$, $0 < \eta \leq 1/L$ and $0 < \beta_2/\beta_1 < 1$.
\end{custheorem}

\begin{cuscorollary}5
\label{coro}
Choose $\eta = \frac{3}{170dL}, b = 25$ and $m = \frac{190dL}{\gamma}$, we have $\frac{\beta_2}{\beta_1} = \frac{3}{4}$.
\end{cuscorollary}
\begin{proof}
With the choice above, we have
\begin{equation}
\beta_1 = 2\eta\left(1-\frac{24}{25}*\frac{3}{170}\right) > 2\eta\cdot\frac{4}{5}
\end{equation}
and
\begin{equation}
\begin{aligned}
\frac{\beta_2}{\beta_1} =& \frac{5}{4}\left(\frac{1}{m\eta\gamma} + \frac{24dL\eta}{bm} + 8dL\eta(\frac{3}{25}+2)\right) \\
=& \frac{5}{4}\left(\frac{17}{57}+(\frac{24}{bm}+\frac{24}{25})dL\eta + 16dL\eta\right) \\
<& \frac{5}{4}\left(\frac{17}{57}+dL\eta+16dL\eta\right) \\
<& \frac{5}{4}\left(\frac{3}{10} + \frac{3}{10}\right) = \frac{3}{4}, \\
\end{aligned}
\end{equation}
where the first inequality comes from
\begin{equation}
\frac{24}{bm} = \frac{24\gamma}{25*190dL} < \frac{24}{25*190d} < \frac{1}{25}
\end{equation}
\end{proof}

\begin{proof} of Theorem \ref{Theorem_1}
For convenience, we define
\begin{equation}
\mathbf{g}_k^s = \frac{1}{\eta}(\mathbf{x}_k^s-\mathbf{x}_{k+1}^s) = \frac{1}{\eta}\left(\mathbf{x}_k^s-Prox_{\eta r}\left(\mathbf{x}_k^s-\eta \hat{\mathbf{v}}_k^s\right)\right)
\end{equation}
so that the proximal gradient step can be written as 
\begin{equation}
\mathbf{x}_{k+1}^s = \mathbf{x}_k^s - \eta\mathbf{g}_k^s
\end{equation}
Then we have
\begin{equation}
\label{t1eq1}
\begin{aligned}
\|\mathbf{x}_{k+1}^s-\mathbf{x}^*\|^2 =& \|\mathbf{x}_k^s-\eta\mathbf{g}_k^s-\mathbf{x}^*\|^2 \\
=& \|\mathbf{x}_k^s-\mathbf{x}^*\|^2 - 2\eta\langle \mathbf{g}_k^s, \mathbf{x}_k^s-\mathbf{x}^*\rangle + \eta^2\|\mathbf{g}_k^s\|^2
\end{aligned}
\end{equation}
Applying Lemma \ref{lemma2} with $\mathbf{x} = \mathbf{x}_k^s, \mathbf{v} = \hat{\mathbf{v}}_k^s, \mathbf{x}^+ = \mathbf{x}_{k+1}^s, \mathbf{g} = \mathbf{g}_k^s$ and $\mathbf{y} = \mathbf{x}^*$, we have
\begin{equation}
-\langle \mathbf{g}_k^s, \mathbf{x}_k^s-\mathbf{x}^*\rangle + \frac{\eta}{2}\|\mathbf{g}_k^s\|^2 \leq F(\mathbf{x}^*)-F(\mathbf{x}_{k+1}^s)- \langle \Delta_k^s, \mathbf{x}_{k+1}^s-\mathbf{x}^*\rangle + L\mu^2,
\end{equation}
where $\Delta_k^s = \hat{\mathbf{v}}_k^s - \nabla f_\mu(\mathbf{x}_k^s)$. Combining with (\ref{t1eq1}), we get
\begin{equation}
\label{t1eq2}
\|\mathbf{x}_{k+1}^s-\mathbf{x}^*\|^2 \leq \|\mathbf{x}_k^s-\mathbf{x}^*\|^2 - 2\eta\left[F(\mathbf{x}_{k+1}^s)-F(\mathbf{x}^*)\right] - 2\eta\langle \Delta_k^s, \mathbf{x}_{k+1}^s-\mathbf{x}^*\rangle + 2\eta L\mu^2
\end{equation}
Now we upper bound the quantity $- 2\eta\langle \Delta_k^s, \mathbf{x}_{k+1}^s-\mathbf{x}^*\rangle$. Denote $\overline{\mathbf{x}}_{k+1}^s = Prox_{\eta r}(\mathbf{x}_k^s - \eta\nabla f_\mu(\mathbf{x}_k^s))$. Then we have
\begin{equation}
\begin{aligned}
- 2\eta\langle \Delta_k^s, \mathbf{x}_{k+1}^s-\mathbf{x}^*\rangle =& - 2\eta\langle \Delta_k^s, \mathbf{x}_{k+1}^s-\overline{\mathbf{x}}_{k+1}^s\rangle - 2\eta\langle \Delta_k^s, \overline{\mathbf{x}}_{k+1}^s-\mathbf{x}^*\rangle \\
\leq& 2\eta\|\mathbf{x}_{k+1}^s-\overline{\mathbf{x}}_{k+1}^s\|\|\Delta_k^s\| - 2\eta\langle \Delta_k^s, \overline{\mathbf{x}}_{k+1}^s-\mathbf{x}^*\rangle \\
\leq& 2\eta\|(\mathbf{x}_k^s-\eta\hat{\mathbf{v}}_k^s)-\left(\mathbf{x}_k^s-\eta\nabla f_\mu(\mathbf{x}_k^s)\right)\|\|\Delta_k^s\| - 2\eta\langle \Delta_k^s, \overline{\mathbf{x}}_{k+1}^s-\mathbf{x}^*\rangle \\
=& 2\eta^2\|\Delta_k^s\|^2 - 2\eta\langle \Delta_k^s, \overline{\mathbf{x}}_{k+1}^s-\mathbf{x}^*\rangle, \\
\end{aligned}
\end{equation}
where the first inequality comes from the Cauchy-Schwarz inequality, the second inequality comes from the non-expandness of proximal mapping. Combining with (\ref{t1eq2}), we get
\begin{equation}
\|\mathbf{x}_{k+1}^s-\mathbf{x}^*\|^2 \leq \|\mathbf{x}_k^s-\mathbf{x}^*\|^2 - 2\eta\left[F(\mathbf{x}_{k+1}^s)-F(\mathbf{x}^*)\right] + 2\eta^2\|\Delta_k^s\|^2 - 2\eta\langle \Delta_k^s, \overline{\mathbf{x}}_{k+1}^s-\mathbf{x}^*\rangle + 2\eta L\mu^2
\end{equation}
Take full expectation on both sides, we obtain
\begin{equation}
\begin{aligned}
\mathbb{E}\left[\|\mathbf{x}_{k+1}^s-\mathbf{x}^*\|^2\right] \leq& \mathbb{E}\left[\|\mathbf{x}_k^s-\mathbf{x}^*\|^2 - 2\eta [F(\mathbf{x}_{k+1}^s)-F(\mathbf{x}^*)] + 2\eta^2\|\Delta_k^s\|^2\right] \\
&- 2\eta\mathbb{E}\left[\langle \Delta_k^s, \overline{\mathbf{x}}_{k+1}^s-\mathbf{x}^*\rangle\right] + 2\eta L\mu^2
\end{aligned}
\end{equation}
From the definition of $\mathbf{v}_k^s$ and Lemma \ref{lemma1} we know that $\mathbb{E}\left[\mathbf{v}_k^s\right] = \nabla f_\mu(\mathbf{x}_k^s)$, i.e., $\mathbb{E}\left[\Delta_k^s\right] = 0$. Then we have
\begin{equation}
\mathbb{E}\left[\langle \Delta_k^s, \overline{\mathbf{x}}_{k+1}^s-\mathbf{x}^*\rangle\right] = \langle \mathbb{E}\left[\Delta_k^s\right], \overline{\mathbf{x}}_{k+1}^s-\mathbf{x}^*\rangle = 0
\end{equation}
In addition, we can bound the term $\mathbb{E}\left[\|\Delta_k^s\|^2\right]$ using Lemma \ref{lemma5} to obtain
\begin{equation}
\begin{aligned}
&\mathbb{E}\left[\|\mathbf{x}_{k+1}^s-\mathbf{x}^*\|^2\right] \\
\leq& \|\mathbf{x}_k^s-\mathbf{x}^*\|^2 - 2\eta\mathbb{E}\left[F(\mathbf{x}_{k+1}^s)-F(\mathbf{x}^*)\right] + 2\eta L\mu^2 + 2\left(\frac{3}{b}+1\right)L^2d^2\mu^2\eta^2 \\
&+ \frac{48dL\eta^2}{b}\mathbb{E}\left[F(\mathbf{x}_k^s)-F(\mathbf{x}^*)\right] + 16dL\eta^2\left(\frac{3}{b}+2\right)\mathbb{E}\left[F(\tilde{\mathbf{x}}_{s-1})-F(\mathbf{x}^*)\right] + 16\eta^2d\|\nabla f(\mathbf{x}^*)\|^2
\end{aligned}
\end{equation}
We now consider a fixed stage $s$, so that $\mathbf{x}_0^s = \tilde{\mathbf{x}}_{s-1}$ and $\tilde{\mathbf{x}}_s = \frac{1}{n}\sum_{k=1}^{m}$. By summing the previous inequality over $k = 0, 1, ..., m-1$ and taking full expectation, we obtain
\begin{equation}
\begin{aligned}
&\mathbb{E}\left[\|\mathbf{x}_m^s-\mathbf{x}^*\|^2\right] + 2\eta\mathbb{E}\left[F(\mathbf{x}_m^s)-F(\mathbf{x}^*)\right] + 2\eta\left(1-\frac{24dL\eta}{b}\right)\sum_{k=1}^{m-1}\mathbb{E}\left[F(\mathbf{x}_k^s)-F(\mathbf{x}^*)\right] \\
\leq& \|\mathbf{x}_0^s-\mathbf{x}^*\|^2 + \frac{48dL\eta^2}{b}\left[F(\mathbf{x}_0^s)-F(\mathbf{x}^*)\right] + 16dL\eta^2\left(\frac{3}{b}+2\right)m\left[F(\tilde{\mathbf{x}}_{s-1})-F(\mathbf{x}^*)\right] \\
&+ 2m\eta L\mu^2\left(1+(\frac{3}{b}+1)\eta Ld^2\right) + 16m\eta^2d\|\nabla f(\mathbf{x}^*)\|^2
\end{aligned}
\end{equation}
Note that $2\eta > 2\eta\left(1-\frac{24d\eta}{b}\right)$ and $\mathbf{x}_0^s = \tilde{\mathbf{x}}_{s-1}$, so we have
\begin{equation}
\begin{aligned}
&2m\eta\left(1-\frac{24dL\eta}{b}\right)\mathbb{E}\left[F(\tilde{\mathbf{x}}_s)-F(\mathbf{x}^*)\right] = 2\eta\left(1-\frac{24dL\eta}{b}\right)\sum_{k=1}^{m}\mathbb{E}\left[F(\mathbf{x}_k^s)-F(\mathbf{x}^*)\right] \\
\leq& \|\mathbf{x}_0^s-\mathbf{x}^*\|^2 + \frac{48dL\eta^2}{b}\left[F(\mathbf{x}_0^s)-F(\mathbf{x}^*)\right] + 16dL\eta^2\left(\frac{3}{b}+2\right)m\left[F(\tilde{\mathbf{x}}_{s-1})-F(\mathbf{x}^*)\right] \\
&+ 2m\eta L\mu^2\left(1+(\frac{3}{b}+1)\eta Ld^2\right) + 16m\eta^2d\|\nabla f(\mathbf{x}^*)\|^2 \\
\leq& \left[\frac{2}{\gamma}+\frac{48dL\eta^2}{b} + 16dL\eta^2\left(\frac{3}{b}+2\right)m\right]\left[F(\tilde{\mathbf{x}}_{s-1})-F(\mathbf{x}^*)\right] \\
&+ 2m\eta L\mu^2\left(1+(\frac{3}{b}+1)\eta Ld^2\right) + 16m\eta^2d\|\nabla f(\mathbf{x}^*)\|^2, \\
\end{aligned}
\end{equation}
where the last inequality comes from the strong convexity of $F$.
Dividing both sides by $m$, we have
\begin{equation}
\begin{aligned}
&2\eta\left(1-\frac{24dL\eta}{b}\right)\mathbb{E}\left[F(\tilde{\mathbf{x}}_s)-F(\mathbf{x}^*)\right] \\
\leq& \left[\frac{2}{m\gamma}+\frac{48dL\eta^2}{bm} + 16dL\eta^2\left(\frac{3}{b}+2\right)\right]\left[F(\tilde{\mathbf{x}}_{s-1})-F(\mathbf{x}^*)\right] \\
&+ 2\eta L\mu^2\left(1+(\frac{3}{b}+1)\eta Ld^2\right) + 16\eta^2d\|\nabla f(\mathbf{x}^*)\|^2\\
\end{aligned}
\end{equation}
Denote $\beta_1 = 2\eta\left(1-\frac{24dL\eta}{b}\right), \beta_2 = \frac{2}{m\gamma}+\frac{48dL\eta^2}{bm} + 16dL\eta^2\left(\frac{3}{b}+2\right)$ and $\delta_\mu = \frac{2\eta L\mu^2}{\beta_1-\beta_2}\left(1+(\frac{3}{b}+1)\eta Ld^2\right)$, then we have
\begin{equation}
\mathbb{E}\left[F(\tilde{\mathbf{x}}_s)-F(\mathbf{x}^*) - \delta_\mu\right] \leq \left(\frac{\beta_2}{\beta_1}\right)\left[F(\tilde{\mathbf{x}}_{s-1})-F(\mathbf{x}^*) - \delta_\mu\right] + \frac{16\eta^2d\|\nabla f(\mathbf{x}^*)\|^2}{\beta_1}
\end{equation}
Telescoping the sum with $s = 1, 2, ..., S$, we get
\begin{equation}
\begin{aligned}
&\mathbb{E}\left[F(\tilde{\mathbf{x}}_S)-F(\mathbf{x}^*) - \delta_\mu\right] \\
\leq& \left(\frac{\beta_2}{\beta_1}\right)^S\left[F(\tilde{\mathbf{x}}_{0})-F(\mathbf{x}^*) - \delta_\mu\right] + \frac{16\eta^2d\|\nabla f(\mathbf{x}^*)\|^2}{\beta_1}\left(1+\left(\frac{\beta_2}{\beta_1}\right)+ ...+\left(\frac{\beta_2}{\beta_1}\right)^{S-1}\right) \\
\leq& \left(\frac{\beta_2}{\beta_1}\right)^S\left[F(\tilde{\mathbf{x}}_{0})-F(\mathbf{x}^*) - \delta_\mu\right] + \frac{16\eta^2d\|\nabla f(\mathbf{x}^*)\|^2}{\beta_1-\beta_2}
\end{aligned}
\end{equation}
Then we complete the proof.
\end{proof}

Now we prove Theorem \ref{remark1}.
\begin{custheorem}3
Suppose Assumptions \ref{a61}, \ref{a3} and \ref{a2} hold,  ZOR-ProxSVRG satisfies the ZOOD property with complexity $\mathcal{C}(L, \gamma) = \mathcal{O}\left(n + \frac{dL}{\gamma}\right)$, $\mathcal{K}=\frac{\beta_2}{\beta_1}$, $\mathcal{E} = \mathcal{O}\left(\|\nabla f(\mathbf{x}^*)\|^2\right)$ and \\ $\delta_\mu = \frac{2\eta L\mu^2}{\beta_1-\beta_2}\left(1+(\frac{3}{b}+1)\eta Ld^2\right)$, where $\beta_1 = 2\eta\left(1-\frac{24dL\eta}{b}\right), \beta_2 = \frac{2}{m\gamma}+\frac{48dL\eta^2}{bm} + 16dL\eta^2\left(\frac{3}{b}+2\right)$ and $\mathbf{x}^* = \arg\min_{\mathbf{x}} F(\mathbf{x})$.
\end{custheorem}
\begin{proof}
From Theorem \ref{Theorem_1} and Corollary \ref{coro} we know that ZOR-ProxSVRG satisfies ZOOD property after running for one epoch and the query complexity is $\mathcal{C}(L, \gamma) = \mathcal{O}\left(n + \frac{dL}{\gamma}\right)$, $\mathcal{E} = \mathcal{O}\left(\|\nabla f(\mathbf{x}^*)\|^2\right)$ and
$\delta_\mu = \frac{2\eta L\mu^2}{\beta_1-\beta_2}\left(1+(\frac{3}{b}+1)\eta Ld^2\right)$, where $\beta_1 = 2\eta\left(1-\frac{24dL\eta}{b}\right), \beta_2 = \frac{2}{m\gamma}+\frac{48dL\eta^2}{bm} + 16dL\eta^2\left(\frac{3}{b}+2\right)$. Then we complete the proof.
\end{proof}

\begin{cuscorollary}3
Suppose Assumptions \ref{a61}, \ref{a4} and \ref{a2}  are satisfied, applying \emph{AdaptRdct-C} on ZOR-ProxSVRG and ZOR-ProxSAGA. From Corollary 1, Theorem \ref{remark1} and \ref{remark2}, denote $G_*^2 = \max\left\{G_\emph{C}^2, G_\emph{NC}^2\right\}$, we have that the total FQC to solve problem (\ref{problem}) is $\tilde{\mathcal{O}}\left(n\log\frac{1}{\epsilon}+\frac{d}{\epsilon}\right)$ with  $\epsilon=\mathbb{E}[F(\mathbf{x}_S)-F(\mathbf{x}^*)] -\mathcal{O}\left( \mu^2 + G_*^2\right)$. 
\end{cuscorollary}

\begin{proof}
From Corollary 1 and Theorem~\ref{remark1} we have the total FQC of ZOR-ProxSVRG is 
\begin{align*}
    \sum_{s=0}^{\log \frac{1}{\epsilon}} \mathcal{O} ( n + \frac{d L}{\gamma_s} )  &= \mathcal{O} (n\log \frac{1}{\epsilon} + \frac{dL}{\gamma_0} \sum_{s=0}^{\log \frac{1}{\epsilon}} \frac{1}{\sqrt{\mathcal{K}}^s } ) \\
    & = \mathcal{O} (n\log \frac{1}{\epsilon} + \frac{d}{\epsilon}).
\end{align*}
From Corollary 1 and Theorem~\ref{remark2} we have the total FQC of ZOR-ProxSAGA is 
\begin{align*}
    \sum_{s=0}^{\log \frac{1}{\epsilon}} \tilde{\mathcal{O}} ( n + \frac{d L}{\gamma_s} )  &= \tilde{\mathcal{O}} (n\log \frac{1}{\epsilon} + \frac{dL}{\gamma_0} \sum_{s=0}^{\log \frac{1}{\epsilon}} \frac{1}{\sqrt{\mathcal{K}}^s } ) \\
    & = \tilde{\mathcal{O}} (n\log \frac{1}{\epsilon} + \frac{d}{\epsilon}).
\end{align*}
Above, the notation $\tilde{\mathcal{O}}$ hides the logarithmic dependence on $d$.
\end{proof}

\section*{Appendix D: Proof of Theorem \ref{remark2}}

We first present Lemma \ref{lemma7} and Theorem \ref{Theorem_5} for our ZOR-ProxSAGA.

\begin{cuslemma}2
\label{lemma7}
Suppose Assumptions \ref{a61}, \ref{a4} and \ref{a2} hold. We have
\begin{equation}
\begin{aligned}
&\mathbb{E}\left[\|\hat{\mathbf{v}}^k-\nabla f_\mu(\mathbf{x}^k)\|^2\right] \\
\leq& 8dL\left(\frac{3}{b}+2\right)\left[F(\mathbf{x}^k) - F(\mathbf{x}^*)\right] + \frac{24dL}{bn}\sum_{i=1}^{n}\left[f_i(\mathbf{\phi}_i^k)-f_i(\mathbf{x}^*) - \langle\nabla f_i(\mathbf{x}^*), \mathbf{\phi}_i^k-\mathbf{x}^*\rangle\right] \\
&+ \left(\frac{3}{b}+1\right)L^2d^2\mu^2 + 8d\|\nabla f(\mathbf{x}^*)\|^2
\end{aligned}
\end{equation}
\end{cuslemma}
\begin{proof}
With the definition of $\hat{\mathbf{v}}^k$, we have
\begin{equation}
\begin{aligned}
\hat{\mathbf{v}}^k-\nabla f_\mu(\mathbf{x}^k) =& \frac{1}{b}\sum_{i \in \mathcal{I}_k}\left(\hat{\nabla}f_i(\mathbf{x}^k)-\hat{\nabla}f_i(\mathbf{\phi}_i^k)\right) + \mathbf{v}^k - \nabla f_\mu(\mathbf{x}^k) \\
=& \frac{1}{b}\sum_{i \in \mathcal{I}_k}\left(\hat{\nabla}f_i(\mathbf{x}^k)-\hat{\nabla}f_i(\mathbf{\phi}_i^k)\right) + \mathbf{v}^k - \hat{\nabla}f(\mathbf{x}^k) + \hat{\nabla}f(\mathbf{x}^k) - \nabla f(\mathbf{x}^k) \\
\end{aligned}
\end{equation}
which is
\begin{equation}
\label{l7eq1}
\begin{aligned}
&\mathbb{E}\left[\|\hat{\mathbf{v}}^k-\nabla f_\mu(\mathbf{x}^k)\|^2\right] \\
=& \mathbb{E}\left[\|\frac{1}{b}\sum_{i \in \mathcal{I}_k}\left(\hat{\nabla}f_i(\mathbf{x}^k)-\hat{\nabla}f_i(\mathbf{\phi}_i^k)\right) + \mathbf{v}^k - \hat{\nabla}f(\mathbf{x}^k) + \hat{\nabla}f(\mathbf{x}^k) - \nabla f(\mathbf{x}^k)\|^2\right] \\
\leq& 2\mathbb{E}\left[\|\frac{1}{b}\sum_{i \in \mathcal{I}_k}\left(\hat{\nabla}f_i(\mathbf{x}^k)-\hat{\nabla}f_i(\mathbf{\phi}_i^k)\right) - \left(\hat{\nabla}f(\mathbf{x}^k) - \mathbf{v}^k\right)\|^2\right] + 2\mathbb{E}\left[\|\hat{\nabla}f(\mathbf{x}^k) - \nabla f(\mathbf{x}^k)\|^2\right], \\
\end{aligned}
\end{equation}
where the first inequality comes from $\|\mathbf{a}+\mathbf{b}\|^2 \leq 2\|\mathbf{a}\|^2+2\|\mathbf{b}\|^2$.
Then we have
\begin{equation}
\label{l7eq2}
\begin{aligned}
&\mathbb{E}\left[\|\frac{1}{b}\sum_{i \in \mathcal{I}_k}\left(\hat{\nabla}f_i(\mathbf{x}^k)-\hat{\nabla}f_i(\mathbf{\phi}_i^k)\right) - \left(\hat{\nabla}f(\mathbf{x}^k) - \mathbf{v}^k\right)\|^2\right] \\
\overset{\textrm{\ding{172}}}{\leq}& \frac{1}{bn}\sum_{i=1}^{n}\mathbb{E}\left[\|\hat{\nabla}f_i(\mathbf{x}^k)-\hat{\nabla}f_i(\mathbf{\phi}_i^k)-\left(\hat{\nabla}f(\mathbf{x}^k) - \mathbf{v}^k\right)\|^2\right] \\
\overset{\textrm{\ding{173}}}{\leq}& \frac{1}{bn}\sum_{i=1}^{n}\mathbb{E}\left[\|\hat{\nabla}f_i(\mathbf{x}^k)-\hat{\nabla}f_i(\mathbf{\phi}_i^k)\|^2\right] \\
\overset{\textrm{\ding{174}}}{\leq}& \frac{6d}{bn}\sum_{i=1}^{n}\mathbb{E}\left[\|\nabla f_i(\mathbf{x}^k) - \nabla f_i(\mathbf{x}^*)\|^2 + \|\nabla f_i(\mathbf{\phi}_i^k) - \nabla f_i(\mathbf{x}^*)\|^2\right] + \frac{3L^2d^2\mu^2}{2b} \\
\overset{\textrm{\ding{175}}}{\leq}& \frac{12dL}{b}\left[F(\mathbf{x}^k) - F(\mathbf{x}^*)\right] + \frac{6d}{bn}\sum_{i=1}^{n}\left[\|\nabla f_i(\mathbf{\phi}_i^k) - \nabla f_i(\mathbf{x}^*)\|^2\right] + \frac{3L^2d^2\mu^2}{2b}, \\
\end{aligned}
\end{equation}
where \ding{172} comes from Lemma \ref{lemma4}; \ding{173} comes from $\mathbb{E}\left[\|\mathbf{a}-\mathbb{E}\left[\mathbf{a}\right]\|^2\right] \leq \mathbb{E}\left[\|\mathbf{a}\|^2\right]$; \ding{174} comes from Lemma \ref{lemma3}; \ding{175} comes from Lemma \ref{lemma4}. Given any $i \in \{1, ..., n\}$, consider the function
\begin{equation}
h_i(\mathbf{x}) = f_i(\mathbf{x}) - f_i(\mathbf{x}^*) - \langle\nabla f_i(\mathbf{x}^*), \mathbf{x}-\mathbf{x}^*\rangle
\end{equation}
It is easy to verify $\nabla h_i(\mathbf{x}^*) = 0$ and $\min_\mathbf{x}h_i(\mathbf{x}) = h_i(\mathbf{x}^*) = 0$. Since each $f_i$ is $L$-smooth, $h_i$ is also $L$-smooth. Then we have
\begin{equation}
\frac{\|\nabla h_i(\mathbf{x}) - \nabla h_i(\mathbf{x}^*)\|^2}{2L} \leq h_i(\mathbf{x})-\min_\mathbf{y}h_i(\mathbf{y}) = h_i(\mathbf{x})-h_i(\mathbf{x}^*) = h_i(\mathbf{x}),
\end{equation}
which implies
\begin{equation}
\label{l7eq3}
\|\nabla f_i(\mathbf{x})-\nabla f_i(\mathbf{x}^*)\|^2 \leq 2L\left[f_i(\mathbf{x})-f_i(\mathbf{x}^*) - \langle\nabla f_i(\mathbf{x}^*), \mathbf{x}-\mathbf{x}^*\rangle\right]
\end{equation}
Plugging (\ref{l7eq3}) into (\ref{l7eq2}), we have
\begin{equation}
\label{l7eq4}
\begin{aligned}
&\mathbb{E}\left[\|\frac{1}{b}\sum_{i \in \mathcal{I}_k}\left(\hat{\nabla}f_i(\mathbf{x}^k)-\hat{\nabla}f_i(\mathbf{\phi}_i^k)\right) - \left(\hat{\nabla}f(\mathbf{x}^k) - \mathbf{v}^k\right)\|^2\right] \\
\leq& \frac{12dL}{b}\left[F(\mathbf{x}^k) - F(\mathbf{x}^*)\right] + \frac{12dL}{bn}\sum_{i=1}^{n}\left[f_i(\mathbf{\phi}_i^k)-f_i(\mathbf{x}^*) - \langle\nabla f_i(\mathbf{x}^*), \mathbf{\phi}_i^k-\mathbf{x}^*\rangle\right] + \frac{3L^2d^2\mu^2}{2b} \\
\end{aligned}
\end{equation}
For the second term on the right hand side of \ref{l7eq1}, we have
\begin{equation}
\label{l7eq5}
\begin{aligned}
\mathbb{E}\left[\|\hat{\nabla}f(\mathbf{x}^k) - \nabla f(\mathbf{x}^k)\|^2\right] \overset{\textrm{\ding{172}}}{\leq}& 2d\|\nabla f(\mathbf{x}^k)\|^2 + \frac{L^2d^2\mu^2}{2} \\
=& 2d\|\nabla f(\mathbf{x}^k) - \nabla f(\mathbf{x}^*) + \nabla f(\mathbf{x}^*)\|^2 + \frac{L^2d^2\mu^2}{2} \\
\overset{\textrm{\ding{173}}}{\leq}& 4d\|\nabla f(\mathbf{x}^k) - \nabla f(\mathbf{x}^*)\|^2 + 4d\|\nabla f(\mathbf{x}^*)\|^2 + \frac{L^2d^2\mu^2}{2} \\
\overset{\textrm{\ding{174}}}{\leq}& \frac{4d}{n}\sum_{i=1}^n\|\nabla f_i(\mathbf{x}^k) - \nabla f_i(\mathbf{x}^*)\|^2 + 4d\|\nabla f(\mathbf{x}^*)\|^2 + \frac{L^2d^2\mu^2}{2} \\
\overset{\textrm{\ding{175}}}{\leq}& 8dL\left[F(\mathbf{x}^k) - F(\mathbf{x}^*)\right] + 4d\|\nabla f(\mathbf{x}^*)\|^2 + \frac{L^2d^2\mu^2}{2}, \\
\end{aligned}
\end{equation}
where \ding{172} comes from Lemma \ref{lemma1}; \ding{173} comes from Young's inequality; \ding{174} comes from Jensen's inequality;  \ding{175} comes from Lemma \ref{lemma4}. Plugging (\ref{l7eq4}) and (\ref{l7eq5}) into (\ref{l7eq1}), we get
\begin{equation}
\begin{aligned}
&\mathbb{E}\left[\|\hat{\mathbf{v}}^k-\nabla f_\mu(\mathbf{x}^k)\|^2\right] \\
\leq& 8dL\left(\frac{3}{b}+2\right)\left[F(\mathbf{x}^k) - F(\mathbf{x}^*)\right] + \frac{24dL}{bn}\sum_{i=1}^{n}\left[f_i(\mathbf{\phi}_i^k)-f_i(\mathbf{x}^*) - \langle\nabla f_i(\mathbf{x}^*), \mathbf{\phi}_i^k-\mathbf{x}^*\rangle\right] \\
&+ \left(\frac{3}{b}+1\right)L^2d^2\mu^2 + 8d\|\nabla f(\mathbf{x}^*)\|^2
\end{aligned}
\end{equation}
Now we complete the proof.
\end{proof}

\begin{custheorem}6	\label{Theorem_5}
Suppose Assumptions \ref{a61}, \ref{a3} and \ref{a2} hold. By running Algorithm \ref{algo4}, we have
\begin{equation}
\mathbb{E}\left[F(\mathbf{x}^K)-F(\mathbf{x}^*)\right] \leq \delta_\mu + \frac{1}{2\eta}\left(1-\frac{\eta\gamma}{2}\right)^K\left[\left(\frac{2}{\gamma}+2\eta+c\right)\left[F(\mathbf{x}^0)-F(\mathbf{x}^*)\right] - \delta\right] + \frac{16d\|\nabla f(\mathbf{x}^*)\|^2}{\gamma},
\end{equation}
where $b = 6, \eta = \min\left\{\frac{1}{95dL}, \frac{2}{3n\gamma}\right\}$, $c = \frac{96\eta^2ndL}{b\left(2b-\eta n\gamma\right)}$, $\delta_\mu = \frac{2(b+3)L^2d^2\mu^2}{b\gamma}$ and $\delta = 2\eta\delta_\mu$.
\end{custheorem}
\begin{proof}
Similar to the proof of Theorem \ref{Theorem_1}, define
\begin{equation}
\mathbf{g}^k = \frac{1}{\eta}(\mathbf{x}^k-\mathbf{x}^{k+1}) = \frac{1}{\eta}\left(\mathbf{x}^k-Prox_{\eta r}\left(\mathbf{x}^k-\eta \hat{\mathbf{v}}^k\right)\right)
\end{equation}
so that the proximal gradient step can be written as 
\begin{equation}
\mathbf{x}^{k+1} = \mathbf{x}^k - \eta\mathbf{g}^k
\end{equation}
Then we have
\begin{equation}
\label{t5eq1}
\begin{aligned}
\|\mathbf{x}^{k+1}-\mathbf{x}^*\|^2 =& \|\mathbf{x}^k-\eta\mathbf{g}^k-\mathbf{x}^*\|^2 \\
=& \|\mathbf{x}^k-\mathbf{x}^*\|^2 - 2\eta\langle \mathbf{g}^k, \mathbf{x}^k-\mathbf{x}^*\rangle + \eta^2\|\mathbf{g}^k\|^2
\end{aligned}
\end{equation}
Applying Lemma \ref{lemma2} with $\mathbf{x} = \mathbf{x}^k, \mathbf{v} = \hat{\mathbf{v}}^k, \mathbf{x}^+ = \mathbf{x}^{k+1}, \mathbf{g} = \mathbf{g}^k$ and $\mathbf{y} = \mathbf{x}^*$, we have
\begin{equation}
-\langle \mathbf{g}^k, \mathbf{x}^k-\mathbf{x}^*\rangle + \frac{\eta}{2}\|\mathbf{g}^k\|^2 \leq F(\mathbf{x}^*)-F(\mathbf{x}^{k+1})- \langle \Delta^k, \mathbf{x}^{k+1}-\mathbf{x}^*\rangle + L\mu^2,
\end{equation}
where $\Delta^k = \hat{\mathbf{v}}^k - \nabla f_\mu(\mathbf{x}^k)$. Combining with (\ref{t5eq1}), we get
\begin{equation}
\label{t5eq2}
\|\mathbf{x}^{k+1}-\mathbf{x}^*\|^2 \leq \|\mathbf{x}^k-\mathbf{x}^*\|^2 - 2\eta\left[F(\mathbf{x}^{k+1})-F(\mathbf{x}^*)\right] - 2\eta\langle \Delta^k, \mathbf{x}^{k+1}-\mathbf{x}^*\rangle + 2\eta L\mu^2
\end{equation}
Now we upper bound the quantity $- 2\eta\langle \Delta^k, \mathbf{x}^{k+1}-\mathbf{x}^*\rangle$. Denote $\overline{\mathbf{x}}^{k+1} = Prox_{\eta r}(\mathbf{x}^k - \eta\nabla f_\mu(\mathbf{x}^k))$. Then we have
\begin{equation}
\begin{aligned}
- 2\eta\langle \Delta^k, \mathbf{x}^{k+1}-\mathbf{x}^*\rangle =& - 2\eta\langle \Delta^k, \mathbf{x}^{k+1}-\overline{\mathbf{x}}^{k+1}\rangle - 2\eta\langle \Delta^k, \overline{\mathbf{x}}^{k+1}-\mathbf{x}^*\rangle \\
\leq& 2\eta\|\mathbf{x}^{k+1}-\overline{\mathbf{x}}^{k+1}\|\|\Delta^k\| - 2\eta\langle \Delta^k, \overline{\mathbf{x}}^{k+1}-\mathbf{x}^*\rangle \\
\leq& 2\eta\|(\mathbf{x}^k-\eta\hat{\mathbf{v}}^k)-\left(\mathbf{x}^k-\eta\nabla f_\mu(\mathbf{x}^k)\right)\|\|\Delta^k\| - 2\eta\langle \Delta^k, \overline{\mathbf{x}}^{k+1}-\mathbf{x}^*\rangle \\
=& 2\eta^2\|\Delta^k\|^2 - 2\eta\langle\Delta^k, \overline{\mathbf{x}}^{k+1}-\mathbf{x}^*\rangle, \\
\end{aligned}
\end{equation}
where the first inequality comes from the Cauchy-Schwarz inequality, the second inequality comes from the non-expandness of proximal mapping. Combining with (\ref{t5eq2}), we get
\begin{equation}
\|\mathbf{x}^{k+1}-\mathbf{x}^*\|^2 \leq \|\mathbf{x}^k-\mathbf{x}^*\|^2 - 2\eta\left[F(\mathbf{x}^{k+1})-F(\mathbf{x}^*)\right] + 2\eta^2\|\Delta^k\|^2 - 2\eta\langle \Delta^k, \overline{\mathbf{x}}^{k+1}-\mathbf{x}^*\rangle + 2\eta L\mu^2
\end{equation}
Take full expectation on both sides, we obtain
\begin{equation}
\begin{aligned}
&\mathbb{E}\left[\|\mathbf{x}^{k+1}-\mathbf{x}^*\|^2\right] \\
\leq& \mathbb{E}\left[\|\mathbf{x}^k-\mathbf{x}^*\|^2 - 2\eta [F(\mathbf{x}^{k+1})-F(\mathbf{x}^*)] + 2\eta^2\|\Delta^k\|^2 - 2\eta\langle \Delta^k, \overline{\mathbf{x}}^{k+1}-\mathbf{x}^*\rangle\right] + 2\eta L\mu^2
\end{aligned}
\end{equation}
From the definition of $\mathbf{v}^k$ and Lemma \ref{lemma1} we know that $\mathbb{E}\left[\mathbf{v}^k\right] = \nabla f_\mu(\mathbf{x}^k)$, i.e., $\mathbb{E}\left[\Delta^k\right] = 0$. Then we have
\begin{equation}
\mathbb{E}\left[\langle \Delta^k, \overline{\mathbf{x}}^{k+1}-\mathbf{x}^*\rangle\right] = \langle \mathbb{E}\left[\Delta^k\right], \overline{\mathbf{x}}^{k+1}-\mathbf{x}^*\rangle = 0
\end{equation}
In addition, we can bound the term $\mathbb{E}\left[\|\Delta^k\|^2\right]$ using Lemma \ref{lemma7} to obtain
\begin{equation}
\begin{aligned}
&\mathbb{E}\left[\|\mathbf{x}^{k+1}-\mathbf{x}^*\|^2\right] + 2\eta [F(\mathbf{x}^{k+1})-F(\mathbf{x}^*)] \\
\leq& \|\mathbf{x}^k-\mathbf{x}^*\|^2 + 16\eta^2dL\left(\frac{3}{b}+2\right)\left[F(\mathbf{x}^k)-F(\mathbf{x}^*)\right] \\
&+ \frac{48\eta^2 dL}{b}\sum_{i=1}^{n}\left[f_i(\mathbf{\phi}_i^k)-f_i(\mathbf{x}^*) - \langle\nabla f_i(\mathbf{x}^*), \mathbf{\phi}_i^k-\mathbf{x}^*\rangle\right] + 2\eta^2\left(\frac{3}{b}+1\right)L^2d^2\mu^2 + 16\eta^2d\|\nabla f(\mathbf{x}^*)\|^2 \\
\leq& \left(1-\frac{\eta\gamma}{2}\right)\|\mathbf{x}^k-\mathbf{x}^*\|^2 + \eta\left(1+16\eta dL\left(\frac{3}{b}+2\right)\right)\left[F(\mathbf{x}^k)-F(\mathbf{x}^*)\right] \\
&+ \frac{48\eta^2 dL}{b}\sum_{i=1}^{n}\left[f_i(\mathbf{\phi}_i^k)-f_i(\mathbf{x}^*) - \langle\nabla f_i(\mathbf{x}^*), \mathbf{\phi}_i^k-\mathbf{x}^*\rangle\right] + 2\eta^2\left(\frac{3}{b}+1\right)L^2d^2\mu^2 + 16\eta^2d\|\nabla f(\mathbf{x}^*)\|^2, \\
\end{aligned}
\end{equation}
where the second inequality comes from the strong convexity of $F$. With the choice of $\mathbf{\phi}_i^{k+1}$, denote $p = b/n$, we have
\begin{equation}
\begin{aligned}
&\frac{1}{n}\sum_{i=1}^{n}\left[f_i(\mathbf{\phi}_i^{k+1})-f_i(\mathbf{x}^*) - \langle\nabla f_i(\mathbf{x}^*), \mathbf{\phi}_i^{k+1}-\mathbf{x}^*\rangle\right] \\
=& p\left[F(\mathbf{x}^k)-F(\mathbf{x}^*)\right] + \frac{1-p}{n}\sum_{i=1}^{n}\left[f_i(\mathbf{\phi}_i^k)-f_i(\mathbf{x}^*) - \langle\nabla f_i(\mathbf{x}^*), \mathbf{\phi}_i^k-\mathbf{x}^*\rangle\right] \\
\end{aligned}
\end{equation}
Define the Lyapunov function
\begin{equation}
\begin{aligned}
&T^k = \|\mathbf{x}^k-\mathbf{x}^*\|^2 + 2\eta\left[F(\mathbf{x}^k)-F(\mathbf{x}^*)\right] + \frac{c}{n}\sum_{i=1}^{n}\left[f_i(\mathbf{\phi}_i^k)-f_i(\mathbf{x}^*) - \langle\nabla f_i(\mathbf{x}^*), \mathbf{\phi}_i^k-\mathbf{x}^*\rangle\right]
\end{aligned}
\end{equation}
Then we have
\begin{equation}
\begin{aligned}
\mathbb{E}\left[T^{k+1}\right] \leq& \left(1-\frac{\eta\gamma}{2}\right)\|\mathbf{x}^k-\mathbf{x}^*\|^2 + \left[\eta\left(1+16\eta dL\left(\frac{3}{b}+2\right)\right) + cp\right]\left[F(\mathbf{x}^k)-F(\mathbf{x}^*)\right] \\
&+ \left(\frac{48\eta^2 dL}{b}+c\left(1-p\right)\right)\sum_{i=1}^{n}\left[f_i(\mathbf{\phi}_i^k)-f_i(\mathbf{x}^*) - \langle\nabla f_i(\mathbf{x}^*), \mathbf{\phi}_i^k-\mathbf{x}^*\rangle\right] \\
&+ 2\eta^2\left(\frac{3}{b}+1\right)L^2d^2\mu^2 + 16\eta^2d\|\nabla f(\mathbf{x}^*)\|^2 \\
\end{aligned}
\end{equation}
Choose $b > 6$, $\eta = \min\left\{\frac{1}{95dL}, \frac{2}{3n\gamma}\right\}$ and $c = \frac{96\eta^2ndL}{b\left(2b-\eta n\gamma\right)}$, denote $\delta = \frac{4(b+3)\eta L^2d^2\mu^2}{b\gamma}$, we have
\begin{equation}
\mathbb{E}\left[T^{k+1} - \delta\right] \leq \left(1-\frac{\eta\gamma}{2}\right)\left[T^k - \delta\right] + 16\eta^2d\|\nabla f(\mathbf{x}^*)\|^2
\end{equation}
Telescope the sum for $k = 0, 1, ..., K-1$, and we have
\begin{equation}
\begin{aligned}
&\mathbb{E}\left[T^K\right] \\
\leq& \delta + \left(1-\frac{\eta\gamma}{2}\right)^K\left[T^0 - \delta\right] + 16\eta^2d\|\nabla f(\mathbf{x}^*)\|^2 \sum_{k=0}^{K-1}\left(1-\frac{\eta\gamma}{2}\right)^k \\
\leq& \delta + \left(1-\frac{\eta\gamma}{2}\right)^K\left[\|\mathbf{x}^0-\mathbf{x}^*\|^2 + 2\eta\left[F(\mathbf{x}^0)-F(\mathbf{x}^*)\right]\right] \\
&+ c\left(1-\frac{\eta\gamma}{2}\right)^K\left[\left[f(\mathbf{x}^0)-f(\mathbf{x}^*)-\langle\nabla f(\mathbf{x}^*), \mathbf{x}^0-\mathbf{x}^*\rangle\right] - \delta\right] + \frac{32\eta d\|\nabla f(\mathbf{x}^*)\|^2}{\gamma}, \\
\end{aligned}
\end{equation}
where the second inequality comes from the fact that $\mathbf{\phi}_i^0 = \mathbf{x}^0, \; \forall i \in [n]$. By the optimality of $\mathbf{x}^*$, there exists $\mathbf{\xi}^*$ such that $\nabla f(\mathbf{x}^*) + \mathbf{\xi}^* = 0$. Then we have
\begin{equation}
\begin{aligned}
f(\mathbf{x}^0)-f(\mathbf{x}^*)-\langle\nabla f(\mathbf{x}^*), \mathbf{x}^0-\mathbf{x}^*\rangle =& f(\mathbf{x}^0)-f(\mathbf{x}^*)+\langle\mathbf{\xi}^*, \mathbf{x}^0-\mathbf{x}^*\rangle \\
\leq& f(\mathbf{x}^0)-f(\mathbf{x}^*) + r(\mathbf{x}^0)-r(\mathbf{x}^*) \\
=& F(\mathbf{x}^0)-F(\mathbf{x}^*),
\end{aligned}
\end{equation}
where the inequality comes from the convexity of $r$. Then we get
\begin{equation}
\begin{aligned}
&\mathbb{E}\left[T^K\right] \leq \delta + \left(1-\frac{\eta\gamma}{2}\right)^K\left[\left(\frac{2}{\gamma}+2\eta+c\right)\left[F(\mathbf{x}^0)-F(\mathbf{x}^*)\right] - \delta\right] + \frac{32\eta d\|\nabla f(\mathbf{x}^*)\|^2}{\gamma}, \\
\end{aligned}
\end{equation}
where the inequality also comes from the strong convexity of $F$. Since each $f_i$ is convex, we have
\begin{equation}
f_i(\mathbf{\phi}_i^K)-f_i(\mathbf{x}^*)-\langle\nabla f_i(\mathbf{x}^*), \mathbf{\phi}_i^K-\mathbf{x}^*\rangle \geq 0, \; \forall i \in [n]
\end{equation}
Then we get
\begin{equation}
\begin{aligned}
&2\eta\mathbb{E}\left[F(\mathbf{x}^K)-F(\mathbf{x}^*)\right] \leq \mathbb{E}\left[T^K\right] \\
\leq& \delta + \left(1-\frac{\eta\gamma}{2}\right)^K\left[\left(\frac{2}{\gamma}+2\eta+c\right)\left[F(\mathbf{x}^0)-F(\mathbf{x}^*)\right] - \delta\right] + \frac{32\eta d\|\nabla f(\mathbf{x}^*)\|^2}{\gamma} \\
\end{aligned}
\end{equation}
which is
\begin{equation}
\mathbb{E}\left[F(\mathbf{x}^K)-F(\mathbf{x}^*)\right] \leq \delta_\mu + \frac{1}{2\eta}\left(1-\frac{\eta\gamma}{2}\right)^K\left[\left(\frac{2}{\gamma}+2\eta+c\right)\left[F(\mathbf{x}^0)-F(\mathbf{x}^*)\right] - \delta\right] + \frac{16d\|\nabla f(\mathbf{x}^*)\|^2}{\gamma},
\end{equation}
where $\eta = \min\left\{\frac{1}{95dL}, \frac{2}{3n\gamma}\right\}$, $c = \frac{96\eta^2ndL}{b\left(2b-\eta n\gamma\right)}$, $\delta_\mu = \frac{2(b+3)L^2d^2\mu^2}{b\gamma}$ and $\delta = 2\eta\delta_\mu$. Now we complete the proof.
\end{proof}

Now we prove Theorem \ref{remark2}.
\begin{custheorem}4
Suppose Assumptions \ref{a61}, \ref{a3} and \ref{a2} hold, ZOR-ProxSAGA satisfies ZOOD property with $\mathcal{C}(L, \gamma) = \mathcal{O}\left(n + \frac{dL}{\gamma}\log\left(\max\left\{dL, n\gamma\right\}\right)\right)$, $\mathcal{E} = \mathcal{O}\left(d\|\nabla f(\mathbf{x}^*)\|^2\right)$ and $\delta_\mu = \frac{2(b+3)L^2d^2\mu^2}{b\gamma}$, where $\mathbf{x}^* = \arg\min_{\mathbf{x}} F(\mathbf{x})$.
\end{custheorem}
\begin{proof}
From Theorem \ref{Theorem_5} we know that ZOR-ProxSAGA satisfies ZOOD property after running $\mathcal{O}\left( \frac{dL}{\gamma}\log\left(\frac{1}{\eta}\right)\right)$ epochs, with query complexity $\mathcal{C}(L, \gamma) = \mathcal{O}\left(n + \frac{dL}{\gamma}\log\left(\max\left\{dL, n\gamma\right\}\right)\right)$, $\mathcal{E} = \mathcal{O}\left(d\|\nabla f(\mathbf{x}^*)\|^2\right)$ and $\delta_\mu = \frac{2(b+3)L^2d^2\mu^2}{b\gamma}$. Then we complete the proof.
\end{proof}

\begin{cuscorollary}4
Suppose Assumptions \ref{a61}, \ref{a1} and \ref{a2}  are satisfied, applying \emph{AdaptRdct-NC} on ZOR-ProxSVRG and ZOR-ProxSAGA. From Corollary 2, Theorem \ref{remark1} and \ref{remark2}, denote $G_*^2 = \max\left\{G_\emph{C}^2, G_\emph{NC}^2\right\}$, we have that the total FQC for  solving problem (\ref{problem}) is 
$\tilde{\mathcal{O}}\left(\frac{n+d}{\epsilon^2}\right)$ with $\epsilon^2=\mathbb{E}||\nabla F(\mathbf{x}_\alpha)||^2 - \mathcal{O}\left( \mu^2 + G_*^2\right)$.
\end{cuscorollary}

\begin{proof}
From Corollary 2 and Theorem~\ref{remark1}, we have the total FQC of ZO-ProxSVRG is 
\begin{align*}
\mathcal{O}(\frac{1}{\epsilon^2}) \cdot \mathcal{O}(n + \frac{dL}{\sigma}) = \mathcal{O}(\frac{n + d}{\epsilon^2}).
\end{align*}
From Corollary 2 and Theorem~\ref{remark2}, we have the total FQC of ZO-ProxSAGA is 
\begin{align*}
    \mathcal{O}(\frac{1}{\epsilon^2}) \cdot \tilde{\mathcal{O}}(n + \frac{dL}{\sigma}) = \tilde{\mathcal{O}}(\frac{n + d}{\epsilon^2}).
\end{align*}
\end{proof}

\section*{Appendix E: Auxillary Lemmas}
Here we present some auxillary lemmas.
\begin{cuslemma}3	\label{lemma1}
Suppose Assumptions \ref{a61} and \ref{a4} holds. For the smoothing function $f_\mu$ defined in (\ref{f_mu}), we have
		
1) $f_\mu(\mathbf{x}) $ 
is also L-smooth and convex, and
\[
\nabla f_\mu(\mathbf{x}) = \mathbb{E}_{\mathbf{u}}[\hat{\nabla}f(\mathbf{x})]
\]
		
2) $\forall \mathbf{x} \in R^d $,
\[
|f_\mu(\mathbf{x}) - f(\mathbf{x})| \leq \frac{L\mu^2}{2}
\]
\[
\|\nabla f_\mu(\mathbf{x}) - \nabla f(\mathbf{x})\|^2 \leq \frac{\mu^2L^2d^2}{4}
\]
		
3) $\forall \mathbf{x} \in R^d $,
\[
\mathbb{E}_\mathbf{u}\left[\|\hat\nabla f(\mathbf{x}) - \nabla f_\mu(\mathbf{x})\|^2\right] \leq \mathbb{E}_{\mathbf{u}} \left[\|\hat{\nabla}f(\mathbf{x})\|^2\right] \leq 2d\|\nabla f(\mathbf{x})\|^2 + \frac{L^2d^2\mu^2}{2}
\]
		
\end{cuslemma}
\begin{proof}
See [\cite{liu2018zeroth}, Lemma 1] and [\cite{gao2018information}, Lemma 4.1].
\end{proof}

\begin{cuslemma}4	\label{lemma3}
Suppose Assumptions \ref{a61} and  \ref{a4} holds. For any $\mathbf{x}, \mathbf{y}$, and any $i \in [n]$, we have
\begin{equation}
\begin{aligned}
&\mathbb{E}\left[\|\hat{\nabla}f_{i}(\mathbf{x})-\hat{\nabla}f_{i}(\mathbf{y})\|^2\right] \\
\leq& 3d\mathbb{E}\left[(1+\beta)\|\nabla f_{i}(\mathbf{x})-\nabla f_{i}(\mathbf{x}^*)\|^2+(1+\frac{1}{\beta})\|\nabla f_{i}(\mathbf{y})-\nabla f_{i}(\mathbf{x}^*)\|^2\right] + \frac{3L^2d^2\mu^2}{2}
\end{aligned}
\end{equation}
\end{cuslemma}
\begin{proof}
From [\cite{ji2019improved}, Appendix D, Lemma 5], we have
\begin{equation}
\begin{aligned}
&\mathbb{E}\left[\|\hat{\nabla}f_{i}(\mathbf{x})-\hat{\nabla}f_{k}(\mathbf{y})\|^2\right] \\
\leq& 3d\mathbb{E}\left[\|\nabla f_{i}(\mathbf{x})-\nabla f_{i}(\mathbf{y})\|^2\right] + \frac{3L^2d^2\mu^2}{2} \\
\leq& 3d\mathbb{E}\left[(1+\beta)\|\nabla f_{i}(\mathbf{x})-\nabla f_{i}(\mathbf{x}^*)\|^2+(1+\frac{1}{\beta})\|\nabla f_{i}(\mathbf{y})-\nabla f_{i}(\mathbf{x}^*)\|^2\right] + \frac{3L^2d^2\mu^2}{2}
\end{aligned}
\end{equation}
The last inequality follows from $\|\mathbf{a}+\mathbf{b}\|^2 \leq (1+\beta)\|\mathbf{a}\|^2 + (1+\frac{1}{\beta})\|\mathbf{b}\|^2$. Then we complete the proof.
\end{proof}

\begin{cuslemma}5
\label{lemma2}
Suppose Assumptions \ref{a61}, \ref{a4} and \ref{a2} hold. Denote \[\mathbf{x}^+ = Prox_{\eta r}\left(\mathbf{x}-\eta\mathbf{v}\right)\] \[\mathbf{g} = \frac{\mathbf{x}^+-\mathbf{x}}{\eta}\] \[\Delta = \mathbf{v}-\nabla f_\mu(\mathbf{x})\] where $\eta$ is a stepsize satisfying $0 < \eta \leq \frac{1}{L}$. Then we have $\forall \; \mathbf{y} \in \mathbb{R}^d$
\begin{equation}
\begin{aligned}
F(\mathbf{y}) \geq F(\mathbf{x^+}) + \langle\mathbf{g}, \mathbf{y}-\mathbf{x}\rangle + \frac{\eta}{2}\|\mathbf{g}\|^2 + \langle\Delta,\mathbf{x}^+-\mathbf{y}\rangle - L\mu^2
\end{aligned}
\end{equation}
\end{cuslemma}
\begin{proof}
We can write the proximal update $\mathbf{x}^+ = Prox_{\eta r}\left(\mathbf{x} - \eta \mathbf{v}\right)$ more explicitly as
\begin{equation}
\mathbf{x}^+ = \arg\min_{\mathbf{y}}\left\{r(\mathbf{y}) + \frac{1}{2\eta}\|\mathbf{y}-\left(\mathbf{x}-\eta\mathbf{v}\right)\|^2\right\}
\end{equation}
The associated optimality condition states that there is a $\mathbf{\xi} \in \partial r(\mathbf{x}^+)$ that
\begin{equation}
\eta\mathbf{\xi} + \mathbf{x}^+ - \left(\mathbf{x}-\eta\mathbf{v}\right) = 0
\end{equation}
Combining with the definition of $\mathbf{g} = \frac{\mathbf{x}-\mathbf{x}^+}{\eta}$, we have $\mathbf{\xi} = \mathbf{g} - \mathbf{v}$.
With the definition of $F$, we have
\begin{equation}
\label{l2eq1}
\begin{aligned}
F(\mathbf{y}) &= f(\mathbf{y}) + r(\mathbf{y}) \\
&\geq f_\mu(\mathbf{y}) - \frac{L\mu^2}{2} + r(\mathbf{y}) \\
&\geq f_\mu(\mathbf{x}) + \langle\nabla f_\mu(\mathbf{x}), \mathbf{y}-\mathbf{x}\rangle + r(\mathbf{x}^+) + \langle\mathbf{\xi}, \mathbf{y}-\mathbf{x}^+\rangle - \frac{L\mu^2}{2}
\end{aligned}
\end{equation}
where the first inequality comes from Lemma \ref{lemma1} and the second inequality comes from the convexity of $f$ and $r$.

From Lemma \ref{lemma1}, we know that $f_\mu$ is also $L$-smooth. Thus we have
\begin{equation}
f_\mu(\mathbf{x}) \geq f_\mu(\mathbf{x}^+) - \langle\nabla f_\mu(\mathbf{x}), \mathbf{x}^+-\mathbf{x}\rangle - \frac{L}{2}\|\mathbf{x}^+-\mathbf{x}\|^2
\end{equation}
Combined with (\ref{l2eq1}), we get
\begin{equation}
\label{l2eq2}
\begin{aligned}
&F(\mathbf{y}) \\
\geq& f_\mu(\mathbf{x}^+) -  \langle\nabla f_\mu(\mathbf{x}), \mathbf{x}^+-\mathbf{x}\rangle - \frac{L}{2}\|\mathbf{x}^+-\mathbf{x}\|^2 + \langle\nabla f_\mu(\mathbf{x}), \mathbf{y}-\mathbf{x}\rangle + r(\mathbf{x}^+) + \langle\mathbf{\xi}, \mathbf{y}-\mathbf{x}^+\rangle - \frac{L\mu^2}{2} \\
\geq& f(\mathbf{x}^+) -  \langle\nabla f_\mu(\mathbf{x}), \mathbf{x}^+-\mathbf{x}\rangle - \frac{L\eta^2}{2}\|\mathbf{g}\|^2 + \langle\nabla f_\mu(\mathbf{x}), \mathbf{y}-\mathbf{x}\rangle + r(\mathbf{x}^+) + \langle\mathbf{\xi}, \mathbf{y}-\mathbf{x}^+\rangle - L\mu^2 \\
\geq& F(\mathbf{x}^+) -  \langle\nabla f_\mu(\mathbf{x}), \mathbf{x}^+-\mathbf{x}\rangle - \frac{L\eta^2}{2}\|\mathbf{g}\|^2 + \langle\nabla f_\mu(\mathbf{x}), \mathbf{y}-\mathbf{x}\rangle + \langle\mathbf{\xi}, \mathbf{y}-\mathbf{x}^+\rangle - L\mu^2 \\
\end{aligned}
\end{equation}
where the second inequality comes from Lemma \ref{lemma1} and the definition that $\mathbf{x}-\mathbf{x}^+=\eta\mathbf{g}$, the third inequality comes from the definition of $F$.
Collecting all inner products on the right-hand side, we have
\begin{equation}
\label{l2eq3}
\begin{aligned}
&- \langle\nabla f_\mu(\mathbf{x}), \mathbf{x}^+-\mathbf{x}\rangle + \langle\nabla f_\mu(\mathbf{x}), \mathbf{y}-\mathbf{x}\rangle + \langle\mathbf{\xi}, \mathbf{y}-\mathbf{x}^+\rangle \\
=& \langle\nabla f_\mu(\mathbf{x}), \mathbf{y}-\mathbf{x}^+\rangle + \langle\mathbf{\xi}, \mathbf{y}-\mathbf{x}^+\rangle \\
\overset{\textrm{\ding{172}}}{=}& \langle\mathbf{g}, \mathbf{y}-\mathbf{x}^+\rangle + \langle\mathbf{v}-\nabla f_\mu(\mathbf{x}), \mathbf{x}^+-\mathbf{y}\rangle \\
\overset{\textrm{\ding{173}}}{=}& \langle\mathbf{g}, \mathbf{y}-\mathbf{x}+\mathbf{x}-\mathbf{x}^+\rangle + \langle\Delta, \mathbf{x}^+-\mathbf{y}\rangle \\
\overset{\textrm{\ding{174}}}{=}& \langle\mathbf{g}, \mathbf{y}-\mathbf{x}\rangle + \eta\|\mathbf{g}\|^2 + \langle\Delta, \mathbf{x}^+-\mathbf{y}\rangle
\end{aligned}
\end{equation}
Above, \ding{172} comes from the definition that $\mathbf{\xi}=\mathbf{g}-\mathbf{v}$, \ding{173} uses $\Delta=\mathbf{v}-\nabla f_\mu(\mathbf{x})$, and \ding{174} uses $\mathbf{x}-\mathbf{x}^+=\eta\mathbf{g}$. Combining (\ref{l2eq2}) and (\ref{l2eq3}), we have
\begin{equation}
\begin{aligned}
F(\mathbf{y}) \geq& F(\mathbf{x}^+) + \langle\mathbf{g}, \mathbf{y}-\mathbf{x}\rangle + \frac{\eta}{2}(2-L\eta)\|\mathbf{g}\|^2 + \langle\Delta, \mathbf{x}^+-\mathbf{y}\rangle - L\mu^2 \\
\geq& F(\mathbf{x}^+) + \langle\mathbf{g}, \mathbf{y}-\mathbf{x}\rangle + \frac{\eta}{2}\|\mathbf{g}\|^2 + \langle\Delta, \mathbf{x}^+-\mathbf{y}\rangle - L\mu^2 \\
\end{aligned}
\end{equation}
where the last inequality holds since we assume $0 < \eta \leq \frac{1}{L}$. Then we complete the proof.
\end{proof}

\begin{cuslemma}6
\label{lemma6}
Let $\left\{\mathbf{z}_i\right\}_{i=1}^n$ be an arbitrary sequence of n vectors. Let $\mathcal{I}$ be a mini-batch of size $b$, which contains $i.i.d.$ samples selected uniformly randomly (with replacement) from $[n]$. Then we have
\begin{equation}
\mathbb{E}_\mathcal{I}\left[\frac{1}{b}\sum_{i \in \mathcal{I}}\mathbf{z}_i\right] = \mathbb{E}\left[\mathbf{z}_i\right] = \frac{1}{n}\sum_{j=1}^{n}\mathbf{z}_j
\end{equation}
When $\mathbb{E}\left[\mathbf{z}_i\right] = 0$, we have
\begin{equation}
\mathbb{E}\left[\left\|\frac{1}{b}\sum_{i \in \mathcal{I}}\mathbf{z}_i\right\|^2\right] = \frac{1}{b}\mathbb{E}\left[\|\mathbf{z}_i\|^2\right] = \frac{1}{bn}\sum_{i=1}^{n}\|\mathbf{z}_i\|^2
\end{equation}
\end{cuslemma}
\begin{proof}
See [\cite{liu2018zeroth}, Lemma 4].
\end{proof}
	
\begin{cuslemma}7
\label{lemma4}
Suppose Assumptions \ref{a61}, \ref{a4} and \ref{a2} hold, consider $F(\mathbf{x})$ as defined in (\ref{problem}), and let $\mathbf{x}^* = \arg\min_{\mathbf{x}} F(\mathbf{x})$. Then we have
\begin{equation}
\frac{1}{n}\sum_{i = 1}^{n}\left\|\nabla f_i(\mathbf{x}) - \nabla f_i(\mathbf{x}^*)\right\|^2 \leq 2L\left[F(\mathbf{x})-F(\mathbf{x}^*)\right]
\end{equation}
\end{cuslemma}
\begin{proof}
See [\cite{xiao2014proximal}, Lemma 1].
\end{proof}

\newpage
\bibliography{ref}
\bibliographystyle{apalike}

\end{document}